\documentclass[reqno,11pt]{amsart}
\usepackage{amsmath,amssymb,amsthm,graphicx,a4wide,enumerate,url}
\usepackage[small,bf]{caption} 
\usepackage{subcaption}
\usepackage{amsmath}
\usepackage{color}
\usepackage{lipsum}
\usepackage{amsfonts}
\usepackage{graphicx}
\usepackage{epstopdf}
\usepackage{algorithmic}
\usepackage[outercaption]{sidecap}
\usepackage{amsopn}
\usepackage{enumerate}
\usepackage{array}
\usepackage{makecell}
\usepackage{adjustbox}
\usepackage{calc}
\usepackage{accents}
\usepackage{hyperref}
\usepackage{cleveref}
\usepackage{romannum}
\usepackage{float}
\usepackage{multirow}
\usepackage{algorithm2e}
\usepackage{array, boldline, makecell, booktabs}

\usepackage{xcolor}
\hypersetup{
  colorlinks   = true, 
  urlcolor     = blue, 
  linkcolor    = blue, 
  citecolor   = red 
}

\DeclareMathOperator{\sech}{sech}

\AtBeginDocument{\pagenumbering{arabic}}
\begin{document}
\parskip.9ex

\title[Multiple-Relaxation ImEx Methods for NLS]
{Accurate Solution of the Nonlinear Schr\"{o}dinger Equation via Conservative Multiple-Relaxation ImEx Methods}

\author[A. Biswas]{Abhijit Biswas}
\address[Abhijit Biswas]
{Computer, Electrical, and Mathematical Sciences \& Engineering Division \\ 
King Abdullah University of Science and Technology \\ Thuwal 23955 \\ Saudi Arabia} 
\email{abhijit.biswas@kaust.edu.sa}

\author[D. I. Ketcheson]{David I. Ketcheson}
\address[David I. Ketcheson]
{Computer, Electrical, and Mathematical Sciences \& Engineering Division \\ 
King Abdullah University of Science and Technology \\ Thuwal 23955 \\ Saudi Arabia} 
\email{david.ketcheson@kaust.edu.sa}
\urladdr{https://www.davidketcheson.info}


\begin{abstract}
    The nonlinear Schr\"{o}dinger (NLS) equation possesses an infinite hierarchy of 
    conserved densities and the numerical preservation of some of these quantities is 
    critical for accurate long-time simulations, particularly for multi-soliton 
    solutions. We propose an essentially explicit discretization that conserves one or 
    two of these conserved quantities by combining higher-order Implicit-Explicit (ImEx) 
    Runge-Kutta time integrators with the relaxation technique and adaptive step size control. 
    We show through numerical tests that our mass-conserving method is much more efficient and accurate than
    the widely-used 2nd-order time-splitting pseudospectral approach.  Compared to higher-order
    operator splitting, it gives similar results in general and significantly
    better results near the semi-classical limit.  Furthermore, for some problems
    adaptive time stepping provides a dramatic reduction in cost without sacrificing accuracy.
    We also propose a full discretization that conserves both mass and energy by using a conservative finite element spatial discretization and multiple relaxation in time. 
    Our results suggest that this method provides a qualitative improvement in long-time
    error growth for multi-soliton solutions.
\end{abstract}

\maketitle
\newcommand{\BigO}{{\mathcal{O}}}
\newcommand{\Dt}{\Delta t}
\newcommand{\Real}{\mathbb{R}}
\newcommand{\Complex}{\mathbb{C}}
\newcommand{\kibitz}[2]{\textcolor{#1}{#2}}
\newcommand{\abhi}[1] {\kibitz{red}{[AB says: #1]}}
\newcommand{\david}[1] {\kibitz{blue}{[DK says: #1]}}
\newcommand{\fim}{f}
\newcommand{\fex}{g}
\newcommand{\atol}{\tol_{\rm{abs}}}
\newcommand{\rtol}{\tol_{\rm{rel}}}
\newcommand{\fac}{\alpha}
\newcommand{\err}{\varepsilon}
\newcommand{\tol}{\tau}
\newcommand{\ueps}{u^\epsilon}
\newcommand{\dtnew}{\Delta t^{\rm new}}


\section{Introduction}
\label{Sec:Introduction}
In this work we study the numerical solution of the focusing nonlinear schr\"{o}dinger (NLS) equation,
which appears in many applications such as nonlinear optics, water waves, plasma physics, biomolecular dynamics, etc. and takes the form:
\begin{align}\label{Eq:Herbst_NLS}
i u_t + u_{xx} + \beta |u|^2 u & = 0,
\end{align}
where $u(x,t)$ is a complex-valued function and $\beta>0$.
In the last few decades, this equation has received significant attention both for its mathematical treatment and the challenges it poses for numerical simulation. It governs the behavior of nonlinear waves in dispersive media and has multi-soliton solutions. While this (completely integrable) equation can be solved analytically using the inverse scattering technique \cite{kruskal1967method,miles1981envelope,shabat1972exact} when the initial condition $u(x,0)$ tends to zero as $|x| \to \infty$, the form of the analytical solution is not very explicit except for the special cases of soliton solutions. Moreover, in practice, various modifications of the NLS equation (such as with a driving force)  are of interest and cannot be solved analytically. As a result, it is necessary to develop numerical methods to understand the dynamics of the NLS equation and its relatives. There is an extensive literature treating its numerical solution \cite{antoine2016high,griffiths1984numerical,herbst1985numerical,verwer1986conservative,sanz1983method}.  Special attention is required in order to avoid non-physical blow-up or unstable numerical solutions \cite{herbst1989numerically,mclaughlin1992chaotic}. The NLS \cref{Eq:Herbst_NLS} has an infinite hierarchy of conserved quantities, among which the first two (mass and energy) are:
\begin{subequations} \label{Herbst_NLS_invariants}
    \begin{align}
        N(t)&:=\int_{\Real} |u(x,t)|^2 dx \equiv N(0),\\
        E(t)&:=\int_{\Real} \left( |u_x(x,t)|^2 - \frac{\beta}{2}|u(x,t)|^4\right) dx \equiv E(0) \;.
    \end{align}
\end{subequations}
Studies have shown that numerical conservation of the two invariants \eqref{Herbst_NLS_invariants} can avoid spurious blow-up or numerical chaos \cite{herbst1985numerical,penha1978contemporary}. It is also known that numerical methods that conserve mass and energy produce asymptotically better error growth compared to non-conservative methods \cite{duran2000numerical}. Therefore, it is natural to design numerical methods that conserve some discrete analogues of these two invariant quantities. Several studies based on finite differences, finite elements, or pseudospectral methods have been conducted to design conservative numerical methods for the NLS equation \cite{delfour1981finite,guo1986convergence,herbst1985numerical,tourigny1988investigation}. However, all the proposed schemes are nonlinearly implicit and computationally expensive. A linearly implicit numerical method is proposed in \cite{fei1995numerical}, but this method is not self-starting and requires small time steps to produce accurate numerical solutions. 

Our main interest is in situations where high accuracy and/or long time integration is required.  In this case,
the popular approach of pseudospectral discretization in space is a natural choice.  However, the 1st- or 2nd-order
operator splitting approach usually employed is not effective or efficient for high-accuracy long-time simulations.
Natural alternatives are high-order operator-splitting methods or high-order ImEx Runge-Kutta methods.
In both cases, the combination with Fourier pseudospectral space discretization yields a completely explicit numerical
method.
High-order ImEx integration of the NLS equation has recently been studied in \cite{antoine2016high}, wherein
the methods performed well, with their main drawback being a lack of exact conservation.

In this paper, we apply relaxation in time, in combination with ImEx methods and adaptive time stepping,
in order to conserve mass and/or energy.  We first describe in Section \ref{Sec:Spatial_discretization} the spatial discretizations
we work with: the Fourier pseudospectral method, which conserves mass, and a finite element (FEM) discretization that conserves both mass and energy.

In Section~\ref{Sec:Diff_time_methods} we describe the time integration methods to be used,
including the widely-used operator splitting approach as well as our proposed ImEx Runge-Kutta
relaxation method with adaptivity.
The relaxation approach can be seen as an oblique projection
onto the conservative manifold after each step, and entails solution of a small system of 1-2 nonlinear algebraic
equation once per step \cite{calvo2006preservation,ketcheson2019relaxation,li2023implicit}.
Conservation of both mass and energy is achieved based on the recently-developed multiple-relaxation
technique \cite{calvo2006preservation,Biswas2023}.  We combine this with a traditional
error-based step size control as well as an additional step size control related to the relaxation
itself.

In Section~\ref{sec:pseudospectral-comparison} we compare the newly-proposed approach with a variety of
other methods, including (low- and high-order) operator splitting as well as ImEx RK without relaxation.
We consider two regimes that are of particular
interest.  In the first regime, $\beta$ is chosen so that the solution consists of a small number of interacting solitons.
The second regime is obtained by taking $0 < \epsilon \ll 1$ and rescaling $x \to \frac{x \epsilon}{\sqrt{2}}$ and $t \to t \epsilon$ to obtain
\begin{align}\label{SemiClassical_NLS_Eqn}
    i \epsilon u^{\epsilon}_t + \frac{\epsilon^2}{2}u^{\epsilon}_{xx} + \beta|u^{\epsilon}|^2 u^{\epsilon} & = 0.\;
\end{align}
In physics this is referred to as the semiclassical regime, and 
$\epsilon$ plays the role of the Planck constant. In this regime, the solution tends to exhibit rapid oscillations in space and time, which become increasingly challenging to capture numerically as $\epsilon \to 0$.
It is known that traditional low-order pseudospectral operator splitting methods struggle to
yield an accurate solution when $\beta>0$ and $\epsilon$ is small \cite{bao2003numerical}, so we focus on this
regime.

For multi-soliton problems we find that all high-order time stepping methods significantly outperform
the low-order Strang or Lie-Trotter splitting methods.  Relaxation incurs a modest computational cost
with a corresponding improvement in overall accuracy.  Thus the ImEx RK relaxation approach provides a flexible
and competitive method that is also mass-conservative.
Similar conclusions apply in the semi-classical regime, as long as
$\epsilon$ is not too small.  For very small $\epsilon$, the proposed relaxation approach gives
an increasingly pronounced advantage over the other methods, and adaptive time stepping improves
computational efficiency by as much as two orders of magnitude.

Finally, in Section~\ref{sec:multiple-relaxation-tests} we study the combination of FEM space discretization and multiple
relaxation ImEX time stepping, comprising the first method that conserves both mass and energy
without requiring the solution of a large system of nonlinear equations at each step.
We find that this discretization seems to yield an error that grows linearly in time even for
2- and 3-soliton solutions (whereas solutions with other methods give a quadratically-growing error).
Conclusions and future work are discussed in Section~\ref{Sec:conclusion}.

\section{Spatial Discretization}
\label{Sec:Spatial_discretization}
We consider two spatial discretizations.  The first
is Fourier pseudo-spectral collocation, which is standard and well-known.  The second is a conservative finite element approach (originally introduced in \cite{herbst1985numerical}).

In the time-splitting pseudospectral approach, the NLS equation is semi-discretized as
\begin{align}
     \dot{U}(t)  & = \mathcal{F}^{-1}\left(\frac{-i \epsilon \xi^2}{2} \mathcal{F}(U)\right) + i\beta |U|^2U,
\end{align}
where $\xi$ is the wave number and $\mathcal{F}$ denotes the Fourier transform. 
This system can be solved using ImEx methods or operator-splitting methods, as described below.
Although \emph{implicit} methods are applied to the first term, the implementation is in fact explicit
since the Fourier transform diagonalizes the system.  With operator splitting, the second term can
be integrated exactly.

The finite element semi-discretization can be written as
\begin{align}\label{ODE_system_FEM}
    \dot{U}(t) = -\left[\frac{1}{h^2} \tilde{I}^{-1}SU(t) + \beta F(U)\right] :=f_{\textrm{FEM}}.
\end{align}
Here $U(t)=\left[v_1(t),w_1(t),v_2(t),w_2(t), \cdots, v_m(t),w_m(t)\right]^T$, where $v_j(t)$ and $w_j(t)$ are the real and the imaginary parts of the solution approximation at $x=x_j$ for $j=1,2,\cdots,m$. With the natural boundary conditions as in \cite{herbst1985numerical}, the matrices $\tilde{I}$ and $S$ are given by
\begin{equation}
	\tilde{I} =\begin{pmatrix}
	{\frac{1}{2}I} & {} & {} & {} & {} \\
	{} & {I} & {} & {O} & {} \\
	{} & {} & {\ddots} & {} & {} \\
	{} & {O} & {} & {I} & {} \\
	{} & {} & {} & {} & {\frac{1}{2}I} 
	\end{pmatrix}
	\ \ \text{,} \ \
	S=
	\begin{pmatrix}
	{-A} & {A} & {} & {} & {} \\
	 {A} & {-2A} & {A} & {O} & {}\\
	{} & {\ddots} & {\ddots} & {\ddots} & {}\\
	{} & {O} & {A} & {-2A} & {A}\\
	{} & {}  & {} & {A} & {-A}\\
	\end{pmatrix}\;,
	\end{equation}
	with 
	\begin{equation}
	I=\begin{pmatrix}
	{1}  & {0}  \\
	{O} & {1} \\
	\end{pmatrix}
	\ \ \text{,} \ \
	A=
	\begin{pmatrix}
	{0} & {1} \\
	{-1} & {0}\\
	\end{pmatrix}\;,
\end{equation}
and the nonlinear function F is given by $F = \left[F_1, F_2, \cdots, F_m\right]^T$, where $F_j = U_{j}(t)^TU_{j}(t)AU_{j}(t)$ with $U_j(t) = \left[v_j(t),w_j(t)\right]^T$.
This system can be shown to conserve the discrete analogs of \eqref{Herbst_NLS_invariants} \cite{herbst1985numerical}.  Defining the discrete mass and energy as
\begin{subequations}
\begin{align}
	\eta_1(U(t)) & := \Delta x \sum_j |U_j(t)|^2 \;\\
	\eta_2(U(t)) & :=  \Delta x \sum_j \left(\left|\frac{U_{j+1}(t)-U_j(t)}{\Delta x}\right|^2 - \frac{1}{2}\beta|U_j(t)|^4\right) \;.
\end{align}
\end{subequations}
It can be shown that
    \begin{align}\label{discrete_inv_eqn}
        \frac{d}{dt} \eta_1 = \frac{d}{dt} \eta_2 = 0.
    \end{align}   

\section{Time Discretization}\label{Sec:Diff_time_methods}
We now turn to the question of temporal discretization.  Due to the nature of the NLS equation, it is advantageous to handle the time integration of the two semi-discrete terms differently.
Starting from either of the spatial discretizations above,
let $U(t)$ denote the semi-discrete approximation of $u(x,t)$, let $\fim(U)$ denote the semi-discretization of $i u_{xx}$, and let $\fex(U)_j = i\beta|U_j|^2 U_j$.  Then we can write the semi-discretization in the form
\begin{equation}\label{IVP_1}
         \dot{U}(t)  = \fim\left(U(t)\right)+\fex\left(U(t)\right), \  U(0) = U_{0}; \ U\in \Complex^{m}.
\end{equation}

\subsection{Operator Splitting}
\label{Sec:TSSP_methods}
Perhaps the most common numerical discretization of the NLS equation is the
highly-effective time-splitting pseudo-spectral method. The time-splitting method involves splitting the equation into simpler sub-problems, which are then solved alternately to advance in time. 
We replace the semi-discretization by the two alternating ODEs
\begin{subequations}
    \begin{align}
    \dot{U}(t) & = \fim(U) \;,\label{split_NLS_a} \\
    \dot{U}(t) & = \fex(U) \;, \label{split_NLS_b}
    \end{align}   
\end{subequations}
and denote the exact solution operator for each part by
$e^{\Delta t\fim}$ and $e^{\Delta t\fex}$, respectively. To solve \cref{split_NLS_a}, we discretize the spatial derivative using a Fourier pseudo-spectral method. The resulting ordinary differential equation system can be integrated exactly in the Fourier space. By taking the inverse Fourier transform, we obtain the solution in the physical space. On the other hand, \cref{split_NLS_b} can be solved exactly in time by making use of an important observation that the quantity $|U|^2$ remains unchanged by the \cref{split_NLS_b} \cite{bao2003numerical}. Specifically, we have
\begin{align*}
U_j(t)= U_j(t_n)e^{\left(i \beta |U_j(t_n)|^2(t-t_n)\right)} \;,
\end{align*}
where $t_n$ is the current time step. Using these solution operators, we compute an updated solution by the multiplicative composition
\begin{align}\label{update_rule_SP}
    U^{n+1} =  e^{a_1\Delta t\fim} e^{b_1\Delta t\fex} \cdots e^{a_s\Delta t\fim} e^{b_s\Delta t\fex} U^{n} \;,
\end{align}
where $U^n$ and $U^{n+1}$ are approximate solutions at times $t_n$ and $t_{n+1}$, respectively and $a_k,b_k \in \Real$ for $k = 1,2,\ldots,s$ define a particular time-splitting method. In this work we use the second-order Strang-splitting method and a fourth-order time-splitting method \cite{auzinger2017practical}, whose coefficients are provided in Table~\ref{Table:1_TS_methods_Coeffs}. 
take fractional steps according to a given time-splitting method, combine them, and then advance the solution from $t_n$ to the next time step $t_{n+1}$ using the update rule given by \cref{update_rule_SP}.
The Strang method is standard and perhaps the most commonly-used method in the literature.  We include the 4th-order method in order to compare operator splitting with our high-order ImEx Runge-Kutta based approach (described below) on a more equal footing.

\begin{table}[ht]
    \centering
    \resizebox{13cm}{!}{%
    \begin{tabular}{ | c| c | c | c | c |}
    \hline
        Method & stage & order & $a_i$ & $b_i$\\
        \hlineB{2.5}
        \multirow{2}*{S2} & \multirow{2}*{2} &  \multirow{2}*{2} & $0.5$ & $1$ \\
        \cline{4-5}
                    & & & $0.5$  & $0$  \\
        \hlineB{2.5}
        \multirow{2}*{AK4} & \multirow{5}*{5} &  \multirow{5}*{4} & $0.267171359000977615$  & $-0.361837907604416033$  \\
        \cline{4-5}
                       & & & $-0.0338279096695056672$  & $0.861837907604416033
                      $ \\
        \cline{4-5}
                       & & & $0.5333131013370561044$  & $0.861837907604416033$ \\
        \cline{4-5}
                        & & & $-0.0338279096695056672
                      $  & $-0.361837907604416033$ \\
        \cline{4-5}
                        & & & $0.267171359000977615$  & $0$ \\
        \hline
    \end{tabular}%
    }
    \caption{Coefficients for two Split-Step methods: S2 (second-order Strang splitting method) and AK4 (a fourth-order splitting method).}
    \label{Table:1_TS_methods_Coeffs}
\end{table}

\subsection{Multiple Relaxation ImEx-RK Methods with a Hybrid Adaptive Step Size Control}
We now introduce the novel time discretization for the NLS equation which is the main subject of this paper.
This approach combines multiple existing methodologies, each possessing unique advantages.  We use high-order ImEx Runge-Kutta time integrators, which can be efficiently applied to the NLS equation since the stiff (second-derivative) term is linear.
One of our primary aims is to develop a numerical method that conserves multiple nonlinear invariants for the NLS equation, which is particularly important for multi-soliton solutions. ImEx methods, as shown in \cite{antoine2016high}, do not guarantee the conservation of nonlinear invariants and may produce inaccurate results, especially when dealing with multi-soliton solutions of NLS, which exhibit steep spatial and temporal gradients.  Conserving multiple nonlinear invariants can help to accurately capture the qualitative behavior of multi-soliton solutions, as highlighted in \cite{herbst1985numerical}. Moreover, since multi-soliton solutions change their shape rapidly over time, an adaptive time-stepping strategy is desirable. In conjunction with ImEx time-stepping methods, we adopt the multiple relaxation approach \cite{Biswas2023} for conserving multiple invariants and propose an adaptive hybrid step size control to efficiently deal with rapid changes in the solution.
\subsubsection{ImEx-RK Methods} \label{SubSubSec:ImEx_Methods}
An ImEx-RK method is a special case of additive RK (ARK) methods \cite{kennedy2003additive} that are designed to solve a system with a vector field with several components evolving on different time scales. We make use of ImEx-RK methods that consist of a diagonally Implicit Runge-Kutta (DIRK) method that is applied to the stiff linear operator $f^{\textrm{Im}}$ and an Explicit Runge-Kutta method (ERK) applied to the non-stiff nonlinear term $f^{\textrm{Ex}}$. All the methods we consider have the same abscissa vector ($\vec{c}$) and the weight vector ($\vec{b}$) for both the implicit and explicit parts of the method. Such methods can be represented via the \emph{Butcher tableau}
\begin{equation}\label{B_table_ImEx_RK}
    \renewcommand\arraystretch{1.2}
    \begin{array}
        {c|c|c}
        \vec{c} & A^{\textrm{Im}} & A^{\textrm{Ex}}\\
        \hline
        & \vec{b}^{\,T} & \vec{b}^{\,T}
    \end{array}\;, 
\end{equation}
where the matrix $A^{\textrm{Im}} = (a^{\textrm{Im}}_{ij})\in \Real^{s \times s}$ with $a^{\textrm{Im}}_{ij} = 0$ for $j > i$, $A^{\textrm{Ex}} = (a^{\textrm{Ex}}_{ij})\in \Real^{s \times s}$ with $a^{\textrm{Ex}}_{ij} = 0$ for $j \geq i$ and vector $\vec{b} \in \Real^{s}$.  We assume $\vec{c} = A^{\textrm{Im}}\vec{e} = A^{\textrm{Ex}}\vec{e},$ where $\vec{e}$ is the vector of ones in $\Real^{s}$. Let $U^n$ and $U^{n+1}$ denote the numerical approximations to the true solution at time $t_n$ and $t_{n+1} = t_n+\Delta t$, respectively. Starting from the solution $U^n$ at time $t_n$, an ImEx-RK method computes the approximate solution at time $t_{n+1}$ as 
\begin{subequations}\label{RK_methods}
    \begin{align}
        g_i & = U^n+\Delta t \sum_{j = 1}^{i} a^{\textrm{Im}}_{ij} f^{\textrm{Im}}(g_j) + \Delta t \sum_{j = 1}^{i-1} a^{\textrm{Ex}}_{ij} f^{\textrm{Ex}}(g_j), \ i = 1,2, \ldots,s \;, \label{RK_methods_a}\\
        U^{n+1} & = u^{n}+\Delta t \sum_{j = 1}^{s} b_j \left(f^{\textrm{Im}}(g_j)+f^{\textrm{Ex}}(g_j)\right)\;. \label{RK_methods_b}
    \end{align}
\end{subequations}
Although the method is formally implicit, note that when pseudo-spectral discretization is used in space,
$\fim$ is diagonal and there is no need to solve any system of equations.  For the FEM spatial discretization, the block-diagonal structure of $\fim$ makes it relatively cheap to solve the linear algebraic system.

This work employs two ImEx methods to compute the presented numerical results:
\begin{itemize}
\item ARK3(2)4L[2]SA: A four-stage, third-order additive Runge-Kutta (ARK) method that includes a second-order embedded ARK method \cite[Appendix C]{kennedy2003additive};
\item ARK4(3)6L[2]SA: A six-stage, fourth-order ARK method that includes a third-order embedded ARK method \cite[Appendix C]{kennedy2003additive}.
\end{itemize}
For the sake of brevity, we will refer to them as ImEx3 and ImEx4, respectively, in the remainder of this paper. In this study, we employ a set of shorthand notations to represent various fully discretized methods used for solving the NLS equation. Specifically, we utilize two spatial discretizations: the Fourier pseudo-spectral collocation method, denoted as 'SP', and the conservative finite element discretization, denoted as 'FEM'. To indicate the use of an ImEx method with relaxation or multiple relaxation, we employ 'R' or 'MR', respectively. Additionally, if a method incorporates a hybrid step size error control, we denote it with 'EC'. For instance, the notation FEM-ImEx4(MR)(EC) implies the implementation of a 4th-order ImEx method with multiple relaxation and hybrid error control for the finite element semi-discretization of the NLS equation. These concise notations provide a convenient and standardized way of representing the specific methods employed in our numerical simulations. Implementations for all the numerical examples presented here can be found in \cite{MRRK_NLS_code}.
\subsubsection{Multiple Relaxation ImEx-RK Methods}
\label{SubSubSec:MR_ImEx_Methods}
Although the FEM semi-discretization discussed above conserves mass and energy, the full discretization with ImEx time integration is not conservative.  In order to impose conservation for the full discretization, we employ the multiple relaxation approach \cite{Biswas2023}.
We briefly review the combination of multiple relaxation and ImEx, which was used but not described in full detail in  \cite{Biswas2023}.  A precise formulation of the ImEx method with relaxation to impose conservation of just one invariant can be found in \cite{li2023implicit}. For simplicity, we present the formulation with two nonlinear invariants, as it will be used here. The extension to more than two invariants following the approach in \cite{Biswas2023} is straightforward.

Let us define $G: = (\eta_1, \eta_2)^T$; then the system \eqref{IVP_1} satisfies     \begin{align*}
    \frac{\textrm{d}}{\textrm{dt}}G(U(t)) = \nabla G(U(t))^{T}f(U(t)) = 0 \;.
\end{align*}
We follow the approach described in \cite{Biswas2023} and consider an ImEx method with an embedded pair 
\begin{equation}\label{B_table_ImEx_RK_Embedded}
    \renewcommand\arraystretch{1.2}
    \begin{array}
        {c|c|c}
        \vec{c} & A^{\textrm{Im}} & A^{\textrm{Ex}}\\
        \hline
        & {\vec{b}^{\,1}}^{\,T} & {\vec{b}^{\,1}}^{\,T}\\
        \hline
        & {\vec{b}^{\,2}}^{\,T} & {\vec{b}^{\,2}}^{\,T}
    \end{array}\;, 
\end{equation}
for imposing the conservation of $\eta_1$ and $\eta_2$. We denote these methods by ImEx$(p,q)$, where $p$ is the order of the original method and $q$ is the order of the embedded method. The modified updated solution $U^{n+1}_{\vec{\gamma}}$ is computed as
    \begin{subequations}\label{modified_Update_rule}
        \begin{align}
            U(t_n + (1+\Gamma )\Delta t) \approx U^{n+1}_{\vec{\gamma}}  & := U^{n+1}+\Delta  t \left(  \gamma_{1} d_{1}^{n} +\gamma_{2} d_{2}^{n}\right) \;, \label{modified_Update_rule_a} \\
            d_{k}^{n} & :=\sum_{j = 1}^{s} b^{k}_{j} \left(f^{\textrm{Im}}(g_j)+f^{\textrm{Ex}}(g_j)\right), \ \text{for} \ k = 1,2\;, \label{increments}
        \end{align}
    \end{subequations}
    where $\Gamma = (\gamma_1+\gamma_2)$, $U^{n+1}:=U^{n}+ \Delta t d_{1}^{n}$, $\vec{\gamma} = (\gamma_{1},\gamma_{2}) \in \Real^{2}$, and the stage values $g_i$ are defined in \eqref{RK_methods_a}. We utilize embedded ImEx methods to compute the directions and use the first method defined by ${\vec{b}^{\,1}}^{\,T}$ to compute $U^{n+1}$. Moreover, we denote the solution of the embedded method as $\hat{U}^{n+1} := U^{n} + \Delta t d_{2}^{n}$. To ensure the conservation of the invariant quantities, the relaxation parameters $(\gamma_{1},\gamma_{2})$ are chosen such that the following nonlinear system of equations is satisfied:
    \begin{align} \label{relaxation-equations}
    G(U^{n+1}_{\vec{\gamma}}) = G(U^n).
    \end{align}
    Successfully solving this system for the relaxation parameters $(\gamma_1, \gamma_2)$ at each step yields a new approximation to the solution at time $t^{n+1}_{\vec{\gamma}} = (t_n + (1+\Gamma )\Delta t$.

\subsubsection{Hybrid Adaptive Step Size Control}
\label{SubSubSec:Hybrid_SSC_MR_ImEX}
In \cite{herbst1985numerical}, a variable time-step method is recommended for solving the NLS equation with rapidly changing solutions. So we equip multiple relaxation ImEx-RK methods with step size controllers based on combining two error indicators, which we refer to as hybrid adaptive step size controllers. The standard approach to automatic step size control \cite{hairer1993solving} based on a local error indicator is to compute two approximations of the solution, $U^{n+1}$ and $\hat{U}^{n+1}$, and define a measure of error as:

\begin{equation}\label{eq:loc_err_control}
\err_n=\left(\frac{1}{m} \sum_{i=1}^m\left(\frac{U_i^{n+1}-\widehat{U}_i^{n+1}}{\atol+\rtol \max \left\{\left|U_i^{n+1}\right|,\left|\hat{U}_i^{n+1}\right|\right\}}\right)^2\right)^{1 / 2} \;.
\end{equation}

Here $m$ is the number of spatial grid points, $\atol$ and $\rtol$ are absolute and relative tolerances, which are both set to $10^{-4}$.  The step is accepted if $\err < 1$, and
a new step size is determined via
\begin{align}\label{modified_control_step}
\dtnew = \fac\left(\frac{1}{\err_n}\right)^{\frac{1}{q+1}} \Delta t \;.
\end{align}
Here $\fac<1$ is a safety factor used to avoid frequent step rejections; we set it to $0.9$ in this work.  If $\err_n>1$, the step is rejected and retried with the step size given by \eqref{modified_control_step}. 

We have found that in rare cases the algebraic solver may fail to find a value of $\vec{\gamma}$ that satisfies
\eqref{relaxation-equations}.  To ensure conservation, if the residual of \eqref{relaxation-equations} is
greater than a specified tolerance (taken herein as $10^{-12}$), we retry the current step with a step size half as large.  This condition is checked and applied before the traditional error-based step size control at each step. The complete adaptive time step algorithm (with a final time $T$) is summarized in Algorithm \ref{step-size-control}.

\RestyleAlgo{boxruled}
\begin{algorithm}[H] \label{step-size-control}
\caption{Hybrid adaptive step control algorithm\label{alg}}
\While{$t < T$}
{
    Compute $U^{n+1}$ and $\hat{U}^{n+1}$ by an ImEx method with embedding and a time step $\Delta t$; \\
    Compute $\err_n$ using \eqref{eq:loc_err_control}; \\
    {
    \eIf{ $\err_n < 1$}
        {
        Find the relaxation parameters $\vec{\gamma}$ and compute: \ $\delta^n = ||G(U^{n+1}_{\vec{\gamma}}) - G(U^n)||_2 $;\\
        \eIf{$\delta^n < 10^{-12}$}
            {Accept the step and update the time and solution as: \\
            $t = (t + (1+\Gamma )\Delta t)$ and $U^{n+1}_{\vec{\gamma}}  = U^{n+1}+\Delta  t \left(  \gamma_{1} d_{1}^{n} +\gamma_{2} d_{2}^{n}\right)$ ;\\
            Go to the next iteration with the time step $\Delta t = \dtnew$, given by \eqref{modified_control_step};}
            {Reject the step and retry with the time step $\Delta t = \frac{\Delta t}{2}$\;}
        }
        {Reject the step and retry with the time step $\Delta t = \dtnew$, given by \eqref{modified_control_step}\;}
    }
}
\end{algorithm}


\section{Numerical Comparison of Time Stepping Methods}
\label{sec:pseudospectral-comparison}
In this section we compare three classes of methods.
All schemes use the same Fourier pseudospectral spatial discretization, while employing different
time discretizations:
\begin{enumerate}
    \item Pseudospectral methods with operator splitting (conservative)
    \item Pseudospectral methods with ImEX RK (non-conservative)
    \item Pseudospectral methods with ImEx RK + relaxation (conservative)
\end{enumerate}
The first two classes of methods have been compared in \cite{antoine2016high}, where the main
drawback of the ImEx schemes was their lack of mass conservation.  In principle the ImEx relaxation schemes,
in which relaxation is used to enforce mass conservation,
can provide all the advantages of this class while remedying the drawback.

We compare the numerical methods above on two classes of problems: 
(a) the NLS equation with multi-soliton solutions, and 
(b) the NLS equation in the semiclassical regime. 

\subsection{Multi-soliton solutions}\label{sec:multi-soliton-comparison}
We first consider solutions consisting of a discrete set of solitons.  These are obtained by taking
the NLS equation \cref{Eq:Herbst_NLS} with the initial condition $u(x,0) = \sech x$. For $\beta=2N^2$ with $N=1,2,\cdots$, the solution consists exactly of the bound state of N solitons. Numerical approximation of the single soliton ($\beta=2$) case has been extensively studied; there are several methods available that can accurately approximate the solution over a long period of time \cite{duran2000numerical,griffiths1984numerical,sanz1983method}.  We focus on the more challenging cases of 2- and 3-soliton solutions; for previous work see \cite{antoine2016high,griffiths1984numerical,herbst1985numerical,verwer1986conservative}. 

\begin{figure}[h]
 \centering
 \begin{subfigure}[b]{0.48\textwidth}
     \centering
     \includegraphics[width=1.0\textwidth]{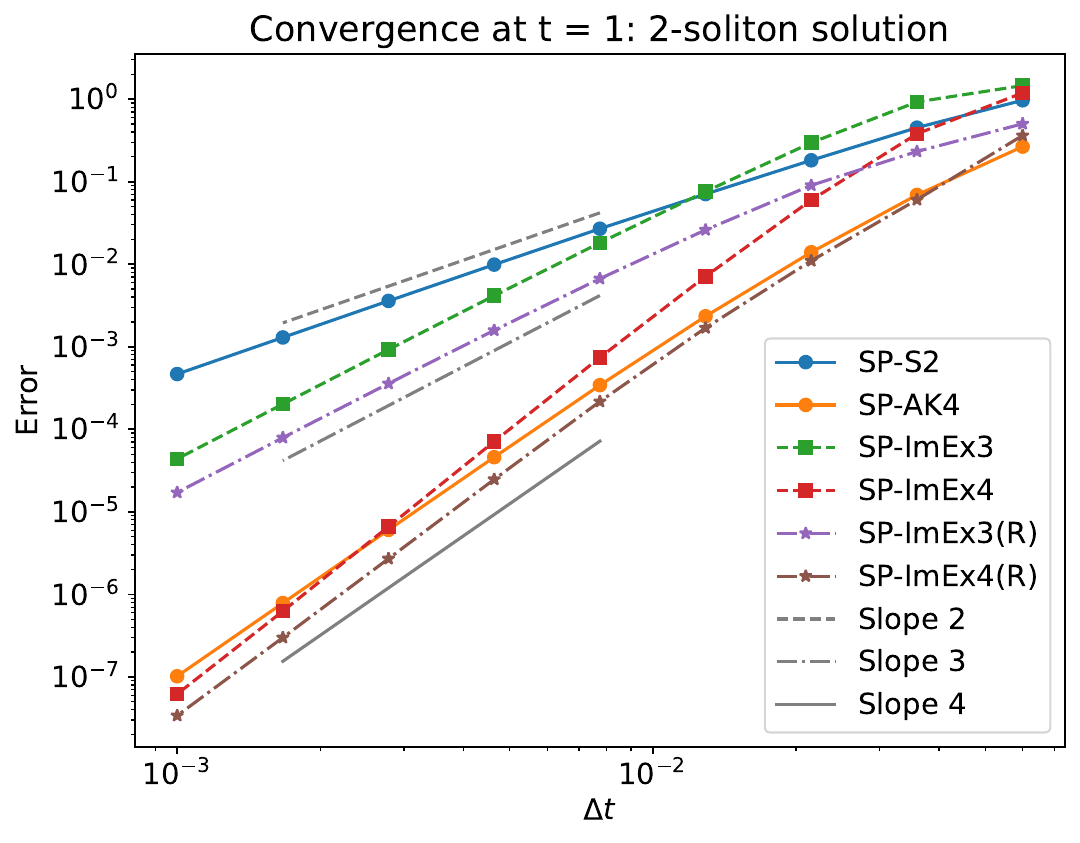}
 \end{subfigure}
 \begin{subfigure}[b]{0.48\textwidth}
     \centering
     \includegraphics[width=1.0\textwidth]{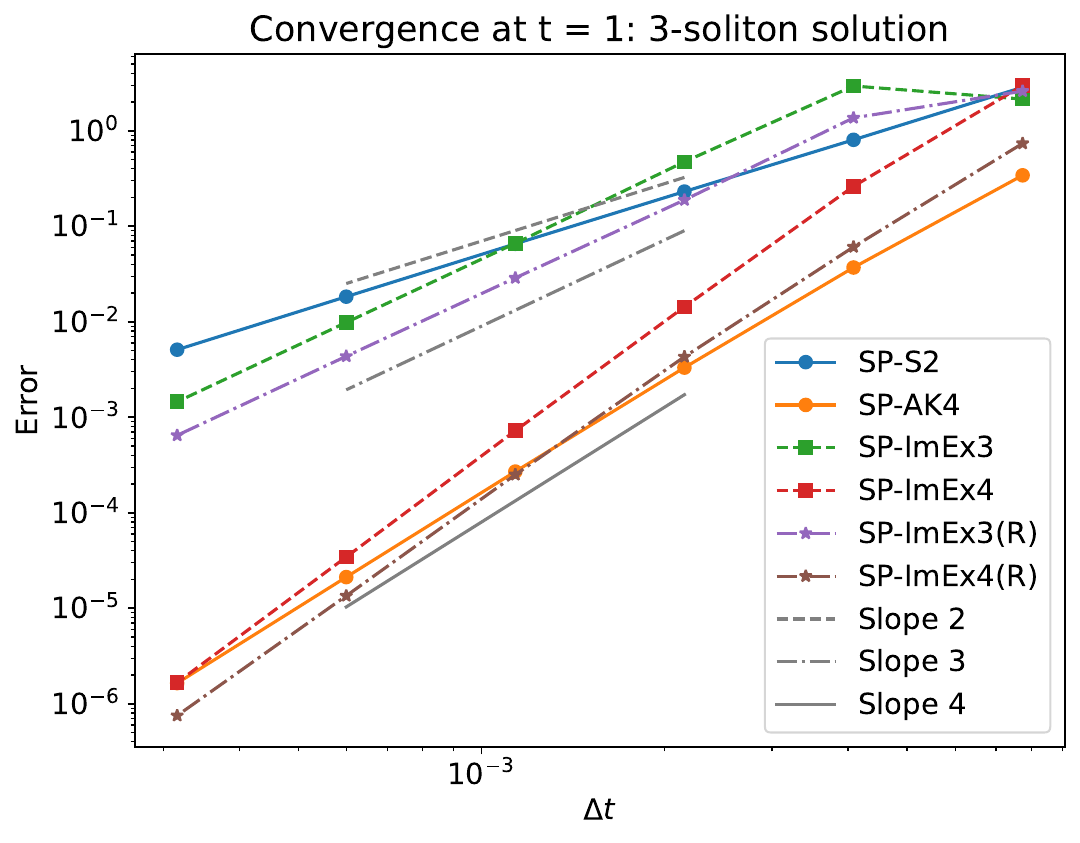}
 \end{subfigure}
	\caption{Convergence of numerical solutions for the NLS \cref{Eq:Herbst_NLS} with 2- and 3-soliton solutions using various time-stepping methods and Fourier pseudo-spectral approximation in space. In the relaxation ImEx methods, conservation of mass is imposed.}
	\label{fig:Sp_space_diff_time_Sol_convg}
\end{figure}

These solutions are not periodic, but by taking a large enough domain we can consider a numerically-equivalent
periodic problem.  To this end,
we take the computational domain $x\in[-35,35]$, which is sufficiently large so that the initial values at the
boundaries are smaller than machine precision.
For the spatial discretization we use $N=1120$ and $N=2240$ grid points, respectively, for the 2- and 3-soliton problems. We solve up to final time $T=1$ and compute the maximum norm of the difference between the numerical and exact solutions at time $T$. To study convergence, we first compute with a sequence of fixed step sizes; results are shown in 
Figure~\ref{fig:Sp_space_diff_time_Sol_convg}. 

These plots confirm that the relaxation approach with ImEx methods retains the same order of accuracy as the original ImEx methods. Furthermore, they demonstrate that applying relaxation to the ImEx methods improves the accuracy of the solutions and yields superior numerical approximations for small $\Delta t$ when compared to time-splitting methods of the same order. Table~\ref{Table:One_Inv_err} presents the errors in mass and energy, and the computational times for various methods, all utilizing a constant starting time step of $\Delta t = 0.01$, and computing the solution up to a final time of $T=5$. Notably, all the relaxation ImEx methods and time-splitting methods exhibit the conservation of mass up to machine precision. Additionally, the computational times for these methods of the same orders are found to be comparable, implying similar computational efficiency.

\begin{table}[ht]
    \centering
    \resizebox{14cm}{!}{%
    \begin{tabular}{ | c| l | c | c | c|}
    \hline
               & &        \multicolumn{2}{c|}{Maximum of} &  \\
   \hline
        Case & Methods & $|\eta_{1}\left(U(t))-\eta_{1}(U(0)\right)|$ & $|\eta_{2}\left(U(t))-\eta_{2}(U(0)\right)|$ & Runtime (s) \\
    \hlineB{2.5}
    \multirow{6}*{2-soliton} & SP-S2  & 6.88E-14 & 1.85E-01  & 1.41 \\
         \cline{2-5}
                             & SP-AK4  & 2.03E-13 & 1.05E-01 & 1.69 \\
         \cline{2-5}                    
                             & SP-ImEx3  & 8.99E-03 & 1.73E-01 & 0.52  \\
         \cline{2-5}
                             & SP-ImEx3(R)  & 1.33E-15 & 9.87E-02 & 0.59 \\
         \cline{2-5}                    
                             & SP-ImEx4  & 5.91E-04 & 1.03E-01 & 0.68  \\
         \cline{2-5}                    
                             & SP-ImEx4(R)  & 1.78E-15 & 1.04E-01 & 0.76 \\         
    \hlineB{2.5}
    \multirow{6}*{3-soliton} & SP-S2  & 1.02E-14 & 5.16E+01  & 1.85 \\
         \cline{2-5}
                             & SP-AK4  & 4.75E-14 & 1.05E+01 & 2.18 \\
         \cline{2-5}                    
                             & SP-ImEx3  & 8.89E-02 & 3.66 & 0.76 \\
         \cline{2-5}
                             & SP-ImEx3(R)  & 8.88E-16 & 2.07 & 0.88 \\
         \cline{2-5}                    
                             & SP-ImEx4  & 8.51E-02 & 5.29 & 1.06 \\
         \cline{2-5}                    
                             & SP-ImEx4(R)  & 4.44E-16 & 1.63 & 1.22 \\ 
    \hlineB{2.5}
    \end{tabular}%
    }
    \caption{Maximum errors in invariants \cref{Herbst_NLS_invariants} by different time integrators applied to the Fourier pseudo-spectral semi-discretization of the NLS equation with different n-solitons. In the relaxation approach with ImEx methods, conservation is exclusively imposed on the mass.}
    \label{Table:One_Inv_err}
\end{table}

\subsection{Semiclassical Regime}
In this section, we consider the problem
\begin{subequations}\label{Semi_Classical_Eq:Herbst_NLS}
    \begin{align}
        i \epsilon u^{\epsilon}_t + \frac{\epsilon^2}{2}u^{\epsilon}_{xx} + |u^{\epsilon}|^2\ u^{\epsilon} & = 0, \text{on} \ [-8,8] \times (0,T] \;, \\
        u^{\epsilon}(-8,t) & =u^{\epsilon}(8,t) \;,
    \end{align}
\end{subequations}
 on the domain $x\in[-8,8]$.  Following \cite{bao2003numerical}, we consider two cases.  First,
 initial data with constant phase: 
 \begin{align} \label{constant-phase}
    \ueps(x,0) & = e^{-x^2};
 \end{align}
and second, initial data with non-uniform phase: 
 \begin{align} \label{varying-phase}
    \ueps(x,0) & = e^{-x^2} e^{i\frac{1}{\epsilon(e^x+e^{-x})}}.
 \end{align}
The same initial conditions were used in \cite{ceniceros2002numerical} to study the semiclassical limits of the NLS equation. In this work, we also consider \cref{Semi_Classical_Eq:Herbst_NLS} with these initial conditions and perform numerical experiments using various time-stepping methods. In the semiclassical regime, one of the important physical quantities of interest is the position density ($\rho^{\epsilon}$), and this is defined by
\begin{align}
   \rho^{\epsilon}(x,t) & = |u^{\epsilon}(x,t)|^2 \;. \label{Eq:SemiclassicalNLS_density} 
\end{align}
Following the meshing strategy for strong $\mathcal{O}(1)$ focusing nonlinearity from \cite{bao2003numerical}, we consider mesh sizes $h = \mathcal{O}(\epsilon)$ and time steps $\Delta t = \small{o}(\epsilon)$, and compare numerical solutions using four methods: the second-order Strang splitting method (S2), the fourth-order splitting method (AK4), the fourth-order Additive Runge-Kutta method (ImEx4), both with and without relaxation. 
For all of these methods we use a fixed time step size.
Additionally, we test the ImEx4 method with both relaxation and error-control (EC) based on step size adaptation,
with a tolerance value of $\tol = 10^{-6}$. In particular, we take the mesh sizes $(\Delta t,\Delta x) = (\frac{1}{100},\frac{1}{32})$ for $\epsilon = 0.2$, $(\Delta t,\Delta x) = (\frac{1}{400},\frac{1}{64})$ for $\epsilon = 0.1$, and $(\Delta t,\Delta x) = (\frac{1}{4000},\frac{1}{128})$ for $\epsilon = 0.05$ and present the numerical approximation for $\rho^{\epsilon}$ (position density) for both initial data in a subinterval at different times. Since an exact solution is not available, to compare our numerical solutions, for each fixed $\epsilon$, we compute the reference "exact" solution by the fourth-order time-splitting method AK4 with a very fine mesh $ (\Delta t,\Delta x) =\left(\frac{1}{10000},\frac{1}{4096}\right)$.

\begin{figure}
     \centering
     \includegraphics[width=\textwidth]{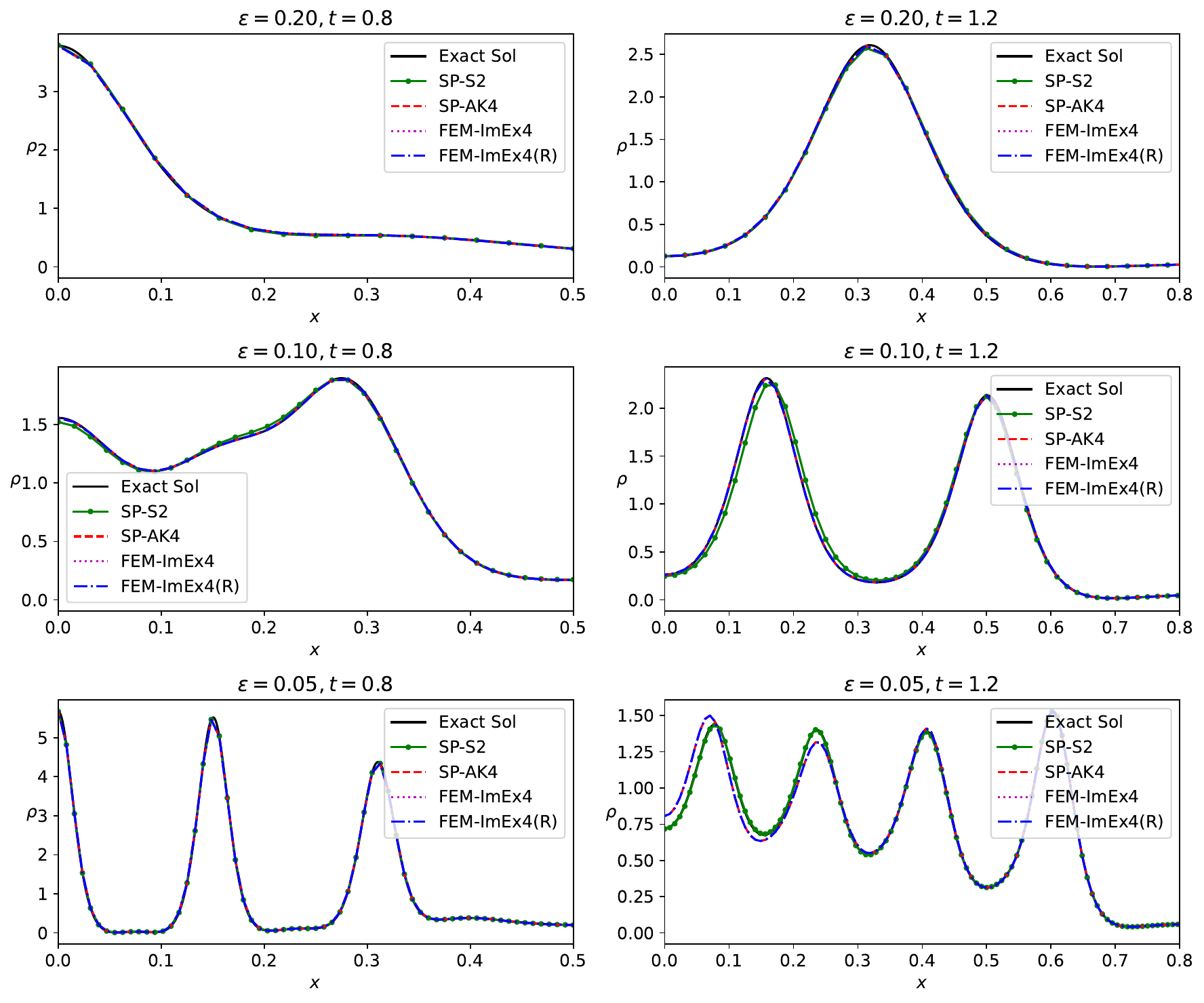}
	\caption{Comparison of numerical density \cref{Eq:SemiclassicalNLS_density} (denoted by $\rho$ in figures) for initial data with constant phase \eqref{constant-phase} and $h = \mathcal{O}(\epsilon)$, $\Delta t = \small{o}(\epsilon)$ meshing strategy. Top row: $\epsilon = 0.2$ and $(\Delta t,\Delta x) = (\frac{1}{100},\frac{1}{32})$; middle row:  $\epsilon = 0.1$ and $(\Delta t,\Delta x) = (\frac{1}{400},\frac{1}{64})$; bottom row: $\epsilon = 0.05$ and $(\Delta t,\Delta x) = (\frac{1}{4000},\frac{1}{128})$.}
   \label{fig:Zero_phase_IC_Semi_NLS}
\end{figure}

In Figure~\ref{fig:Zero_phase_IC_Semi_NLS}, we present the numerical density results for the constant-phase initial data \eqref{constant-phase} at two different times: $t = 0.8$ and $t = 1.2$.
For purposes of comparison we use the same mesh sizes used in \cite{bao2003numerical}, although for the high
order methods used here the errors are generally too small to be seen.  To give more insight,
in Table~\ref{table:Zero_phase_IC_Semi_NLS_err}, we compare the runtimes and the errors in the maximum norm obtained using different methods.  We see that relaxation improves the overall accuracy of the solution, and adaptive error control
dramatically speeds up the computation, particularly when $\epsilon$ is small. It is interesting to note that the SP-S2 method is surprisingly accurate at the smallest value of $\epsilon$.

\begin{table}
    \centering
    \resizebox{15cm}{!}{%
    \begin{tabular}{|l|c|c|c|c|c|c|c|} 
    \hline
        {\begin{tabular}[c]{@{}c@{}} 
        \textbf{Method} \end{tabular}} &
        {\begin{tabular}[c]{@{}c@{}}  
        \textbf{Solution error} \end{tabular}} &
        {\begin{tabular}[c]{@{}c@{}}   \textbf{Runtime (s)} \end{tabular}} &
        {\begin{tabular}[c]{@{}c@{}} 
        \textbf{Solution  error}\end{tabular}} &
        {\begin{tabular}[c]{@{}c@{}} \textbf{Runtime} (s) \end{tabular}}\\
    \hline
    \hline
        &         \multicolumn{2}{c|}{$\epsilon = 0.2$, \ $t= 0.8$} & \multicolumn{2}{c|}{$\epsilon = 0.2$, \ $t= 1.2$}\\
    \hline
       SP-S2 &  7.41E-02	& 0.01 & 6.25E-02 & 0.02 \\
       SP-AK4 &  1.99e-03	& 0.02 & 1.73E-03 & 0.03\\
       FEM-ImEx4 & 9.19E-04	& 0.06 & 1.37E-03 & 0.09 \\
       FEM-ImEx4(R) & 3.49E-04 & 0.07	& 4.13E-04 & 0.11 \\
       FEM-ImEx4(R)(EC) & 2.36E-04 & 0.06	& 1.43E-04 & 0.09 \\
    \hline
            &         \multicolumn{2}{c|}{$\epsilon = 0.1$, \ $t= 0.8$} & \multicolumn{2}{c|}{$\epsilon = 0.1$, \ $t= 1.2$}\\
    \hline

       SP-S2 & 2.68E-02	&  0.09 & 1.32E-01 & 0.17\\
       SP-AK4 & 2.00E-04	&  0.17 & 7.47E-04 & 0.27\\
       FEM-ImEx4 & 4.77E-04 & 0.32 &  1.54E-03 & 0.59 \\
       FEM-ImEx4(R) & 1.85E-04 & 0.37	& 1.31E-03 & 0.65 \\
       FEM-ImEx4(R)(EC) & 2.50E-04	& 0.17 & 1.81E-03 & 0.25 \\
       
    \hline
        &         \multicolumn{2}{c|}{$\epsilon = 0.05$, \ $t= 0.8$} & \multicolumn{2}{c|}{$\epsilon = 0.05$, \ $t= 1.2$}\\   
    \hline
       SP-S2 & 4.40E-03	& 50.63 & 1.29E-02 &  124.63 \\
       SP-AK4 & 1.29E-02	& 53.01 &  2.12E-01 & 129.52\\
       FEM-ImEx4 & 1.29E-02  & 53.02 & 2.12E-01 & 129.92\\
       FEM-ImEx4(R) & 1.29E-02 & 55.60	& 2.12E-01 & 134.48 \\
       FEM-ImEx4(R)(EC) & 6.88E-03 &  0.62 & 1.25E-01 & 1.38\\
    \hline
    \end{tabular}
    }%
    \caption{Error in the maximum norm of the  solution and the runtime in different cases for initial data with constant phase \eqref{constant-phase} and $h = \mathcal{O}(\epsilon)$, $\Delta t = \small{o}(\epsilon)$ meshing strategy. Top row: $\epsilon = 0.2$ and $(\Delta t,\Delta x) = (\frac{1}{100},\frac{1}{32})$; middle row:  $\epsilon = 0.1$ and $(\Delta t,\Delta x) = (\frac{1}{400},\frac{1}{64})$; bottom row: $\epsilon = 0.05$ and $(\Delta t,\Delta x) = (\frac{1}{4000},\frac{1}{128})$.}
    \label{table:Zero_phase_IC_Semi_NLS_err}
\end{table}	

\begin{figure}
 \centering
     \centering
     \includegraphics[width=\textwidth]{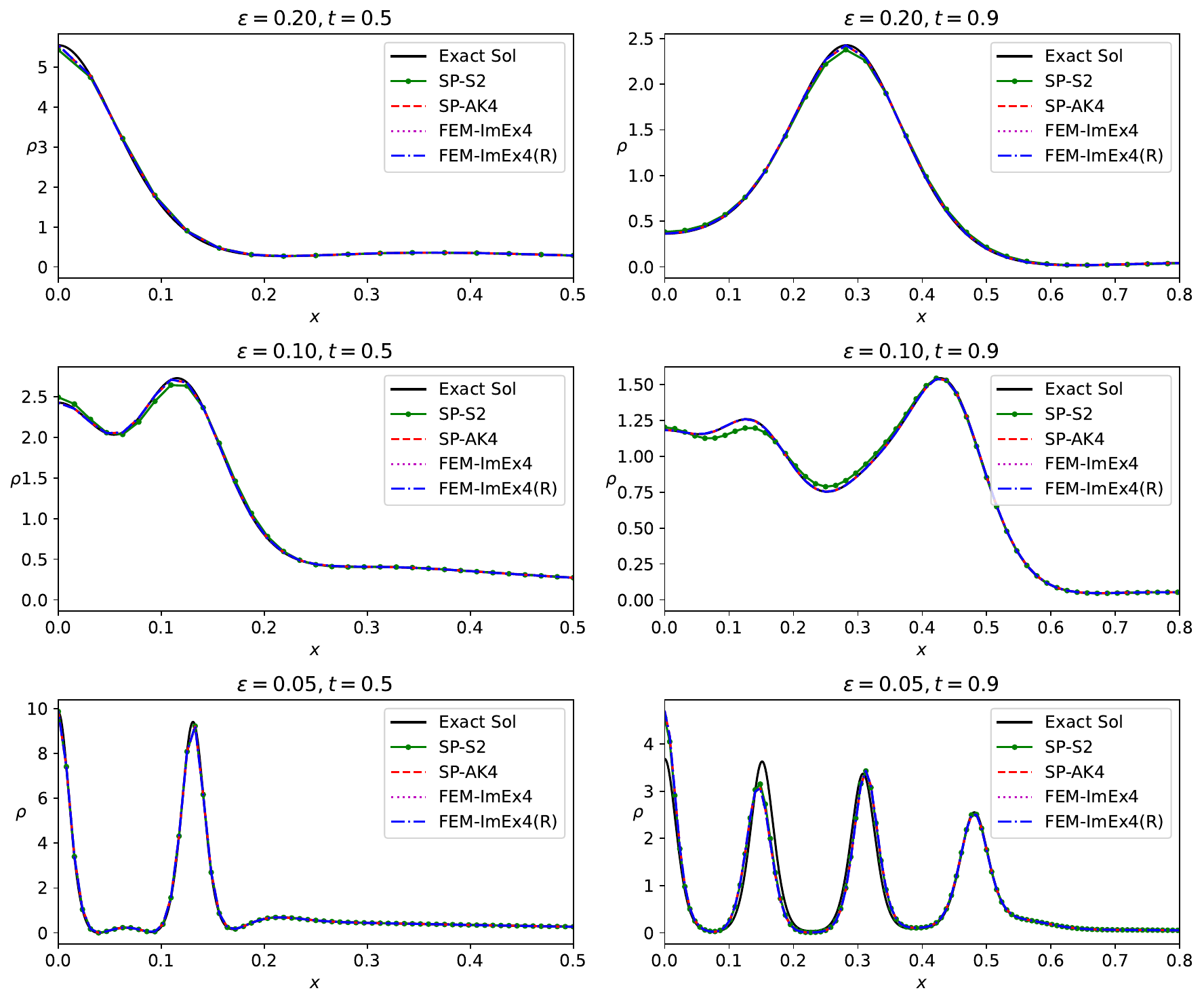}
	\caption{Comparison of numerical density \cref{Eq:SemiclassicalNLS_density} with varying initial phase \eqref{varying-phase} and $h = \mathcal{O}(\epsilon)$, $\Delta t = \small{o}(\epsilon)$ meshing strategy. Top row: $\epsilon = 0.2$ and $(\Delta t,\Delta x) = (\frac{1}{100},\frac{1}{32})$; middle row:  $\epsilon = 0.1$ and $(\Delta t,\Delta x) = (\frac{1}{400},\frac{1}{64})$; bottom row: $\epsilon = 0.05$ and $(\Delta t,\Delta x) = (\frac{1}{4000},\frac{1}{128})$.}
	\label{fig:Nonzero_phase_IC_Semi_NLS}
\end{figure}

\begin{table}
    \centering
    \resizebox{15cm}{!}{%
    \begin{tabular}{|l|c|c|c|c|c|c|c|} 
    \hline
        {\begin{tabular}[c]{@{}c@{}} 
        \textbf{Method}
        \end{tabular}} &
        {\begin{tabular}[c]{@{}c@{}}  
        \textbf{Solution error} \end{tabular}} &
        {\begin{tabular}[c]{@{}c@{}}   \textbf{Runtime (s)} \end{tabular}} &
        {\begin{tabular}[c]{@{}c@{}} 
        \textbf{Solution  error}\end{tabular}} &
        {\begin{tabular}[c]{@{}c@{}} \textbf{Runtime} (s) \end{tabular}}\\
    \hline
    \hline
        &         \multicolumn{2}{c|}{$\epsilon = 0.2$, \ $t= 0.5$} & \multicolumn{2}{c|}{$\epsilon = 0.2$, \ $t= 0.9$}\\ 
    \hline
       SP-S2 &  2.70E-02	& 0.01 & 6.26E-02 & 0.01 \\
       SP-AK4 &  1.25E-03	& 0.01 & 2.11E-03 & 0.03\\
       FEM-ImEx4 & 2.14E-04	& 0.04 & 1.38E-03 & 0.06 \\
       FEM-ImEx4(R) & 1.84E-04 & 0.04	& 4.34E-04 & 0.08 \\
       FEM-ImEx4(R)(EC) & 2.92E-04 & 0.03	& 9.75E-05 & 0.07 \\
    \hline
        &         \multicolumn{2}{c|}{$\epsilon = 0.1$, \ $t= 0.5$} & \multicolumn{2}{c|}{$\epsilon = 0.1$, \ $t= 0.9$}\\ 
        
    \hline
       SP-S2  & 4.41E-02	&  0.05 & 1.06E-01 & 0.11\\
       SP-AK4 & 1.16E-04	&  0.09 & 2.35E-04 & 0.18\\
       FEM-ImEx4 & 6.42E-04 & 0.20 &  9.20E-04 & 0.37 \\
       FEM-ImEx4(R) & 8.62E-04 & 0.23	& 2.36E-03 & 0.44 \\
       FEM-ImEx4(R)(EC) & 1.04E-03	& 0.09 & 3.18E-03 & 0.22 \\
       
    \hline
        &         \multicolumn{2}{c|}{$\epsilon = 0.05$, \ $t= 0.5$} & \multicolumn{2}{c|}{$\epsilon = 0.05$, \ $t= 0.9$}\\  
    \hline
       SP-SP2 & 2.08E-02	& 21.12 & 1.56 &  72.02 \\
       SP-AK4 & 3.40E-02	& 21.51 & 1.92 & 73.05\\
       FEM-ImEx4 & 3.40E-02  & 23.34 & 1.92 & 70.77\\
       FEM-ImEx4(R) & 3.40E-02 & 24.06	& 1.92 & 79.41 \\
       FEM-ImEx4(R)(EC) & 3.04E-02 &  0.32 & 1.86 & 1.04\\
    \hline
    \end{tabular}
    }%
    \caption{Error in the maximum norm of the solution and the runtime in different cases for the initial data with varying initial phase \eqref{varying-phase} and $h = \mathcal{O}(\epsilon)$, $\Delta t = \small{o}(\epsilon)$ meshing strategy. Top row: $\epsilon = 0.2$ and $(\Delta t,\Delta x) = (\frac{1}{100},\frac{1}{32})$; middle row:  $\epsilon = 0.1$ and $(\Delta t,\Delta x) = (\frac{1}{400},\frac{1}{64})$; bottom row: $\epsilon = 0.05$ and $(\Delta t,\Delta x) = (\frac{1}{4000},\frac{1}{128})$.}
    \label{table:Nonzero_phase_IC_Semi_NLS_err}
\end{table}

In Figure~\ref{fig:Nonzero_phase_IC_Semi_NLS}, we present the numerical densities for the nonzero initial phase data at time $t = 0.5$ and $t = 0.9$. Table~\ref{table:Nonzero_phase_IC_Semi_NLS_err} provides a comparison of the runtimes and errors in the maximum norm obtained by different methods for the non-zero initial phase data. In this scenario, all fourth-order methods exhibit similar errors, but again the application of step-size adaptivity yields a huge reduction in computational cost
for small $\epsilon$.

\begin{figure}
    \centering
    \begin{subfigure}[b]{0.48\textwidth}
            \includegraphics[width=\textwidth]{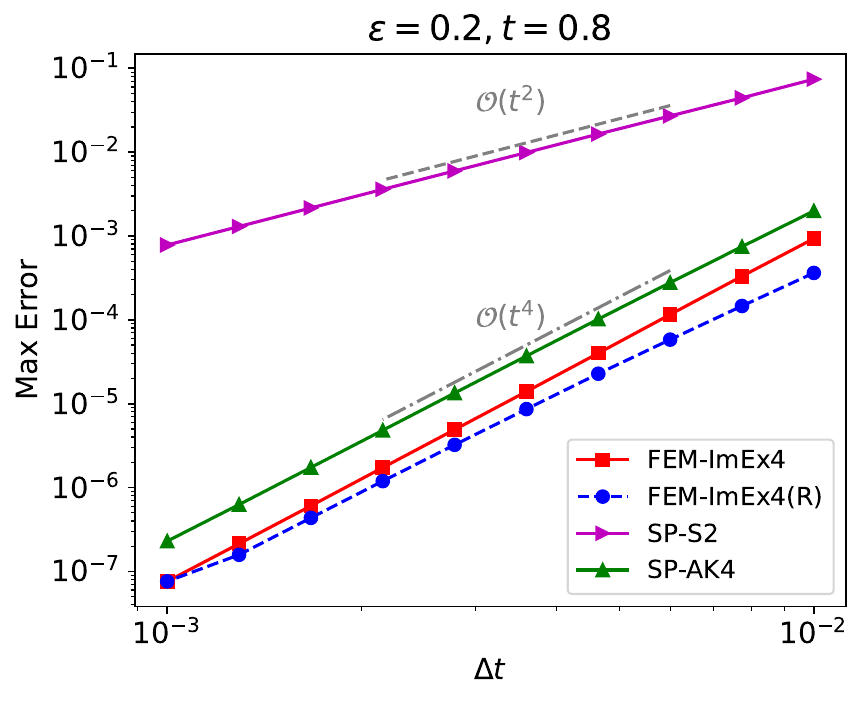}
            \caption{Solution error convergence.}
    \end{subfigure}
    \begin{subfigure}[b]{0.48\textwidth}
            \includegraphics[width=\textwidth]{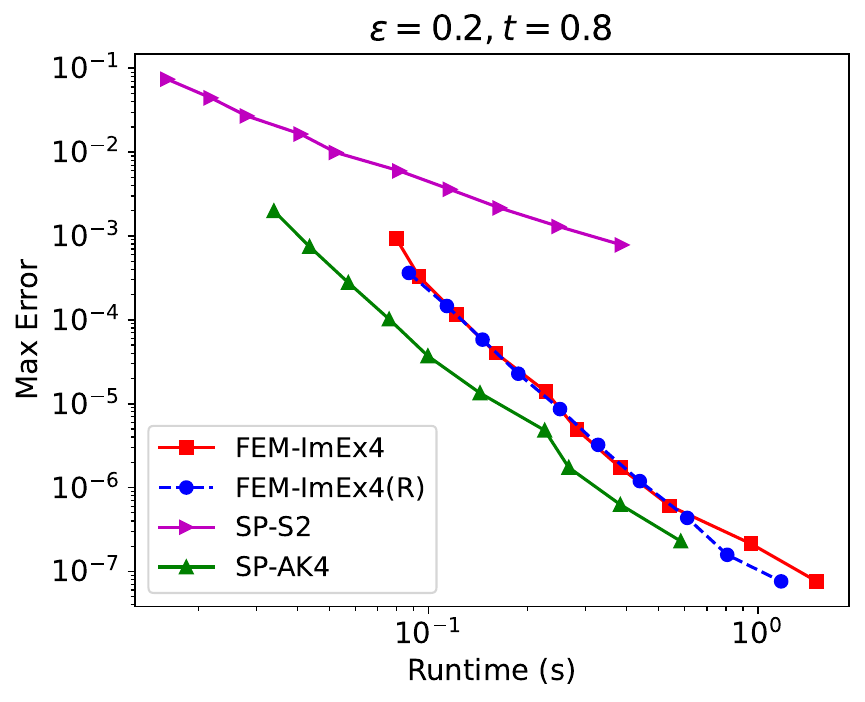}
            \caption{Cost measurement comparison.}
    \end{subfigure}
    \begin{subfigure}[b]{0.48\textwidth}
            \includegraphics[width=\textwidth]{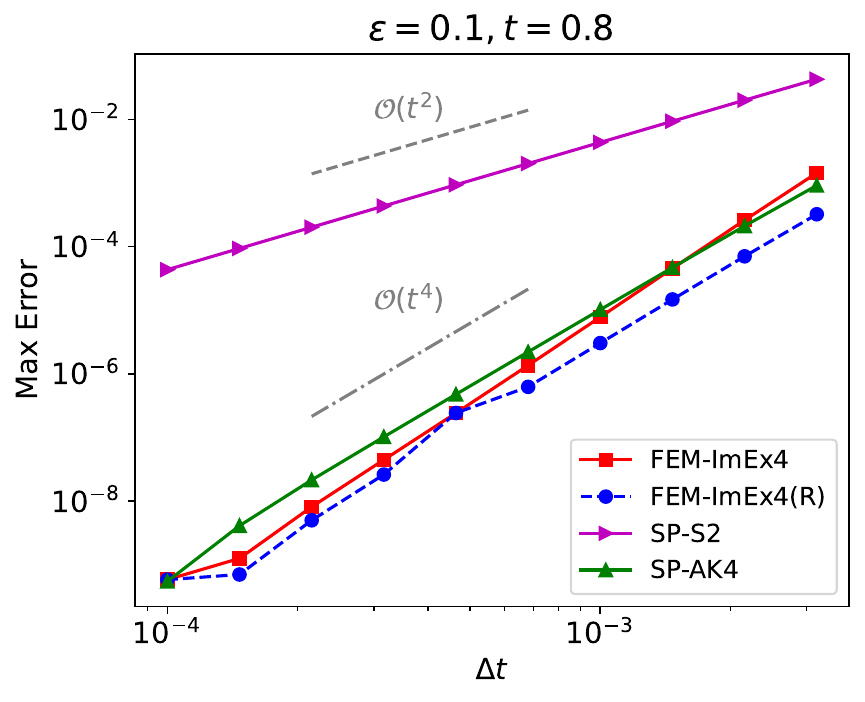}
            \caption{Solution error convergence.}
    \end{subfigure}
    \begin{subfigure}[b]{0.48\textwidth}
            \includegraphics[width=\textwidth]{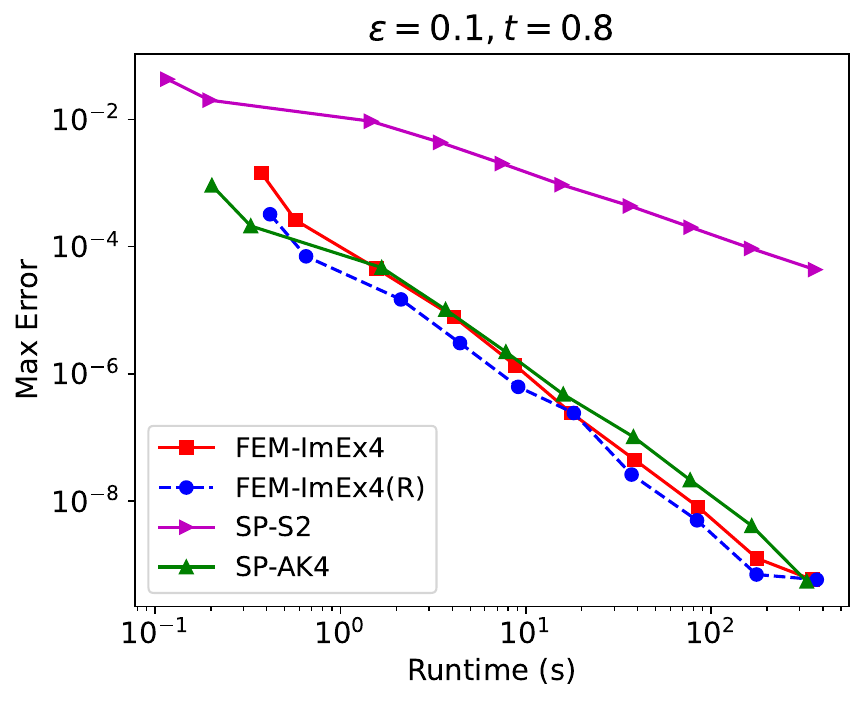}
            \caption{Cost measurement comparison.}
    \end{subfigure}
    \caption{Convergence rates and efficiency of different methods for the strong $\mathcal{O}(1)$ cubic focusing nonlinearity with $\epsilon = 0.2, 0.1$ and a zero initial phase. Meshing strategies for top row: $(\epsilon, \Delta x) = \left(0.2,\frac{1}{64}\right)$; bottom row: $(\epsilon,\Delta x) = \left(0.1,\frac{1}{128}\right)$.}
    \label{fig:Zero_phase_IC_cost_comp}
\end{figure}

\begin{figure}
    \centering
    \begin{subfigure}[b]{0.48\textwidth}
            \includegraphics[width=\textwidth]{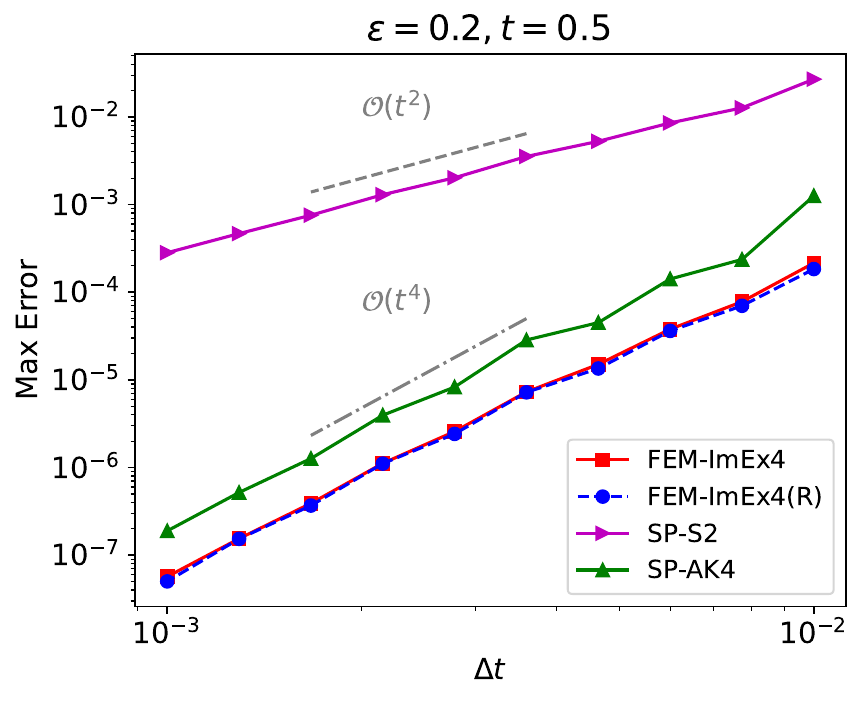}
            \caption{Solution error convergence.}
    \end{subfigure}
    \begin{subfigure}[b]{0.48\textwidth}
            \includegraphics[width=\textwidth]{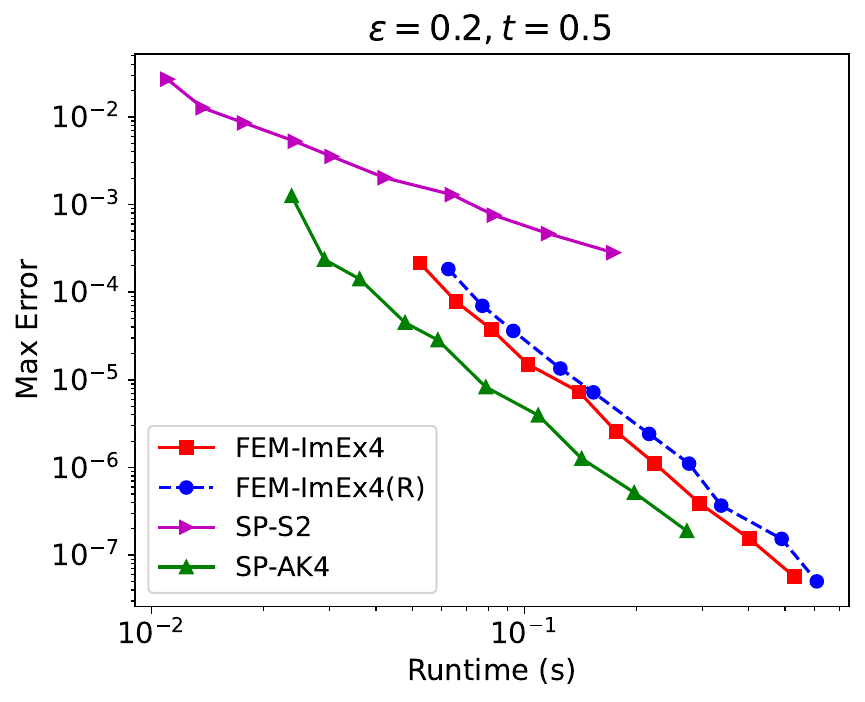}
            \caption{Cost measurement comparison.}
    \end{subfigure}
    \begin{subfigure}[b]{0.48\textwidth}
            \includegraphics[width=\textwidth]{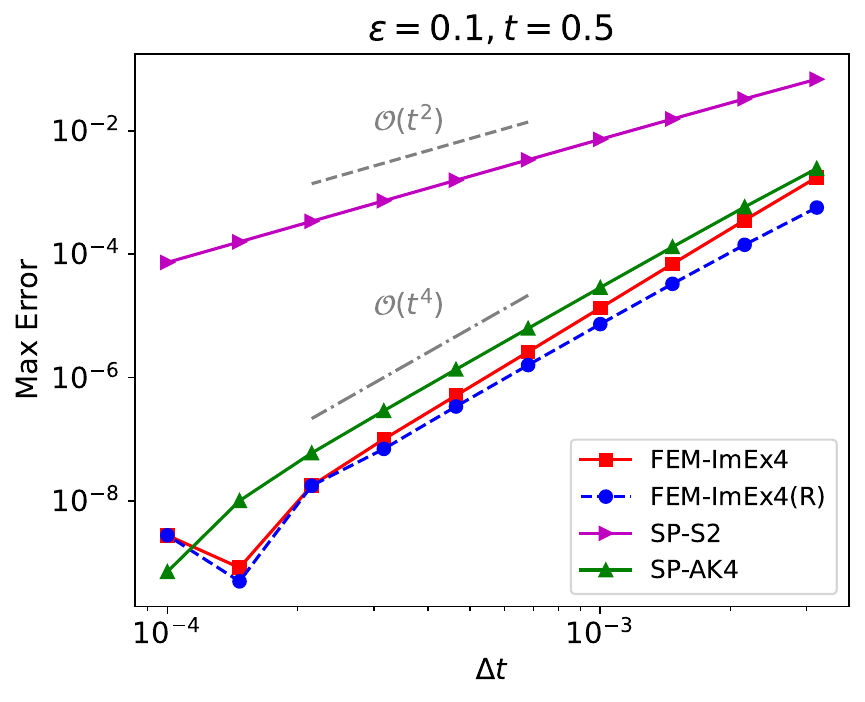}
            \caption{Solution error convergence.}
    \end{subfigure}
    \begin{subfigure}[b]{0.48\textwidth}
            \includegraphics[width=\textwidth]{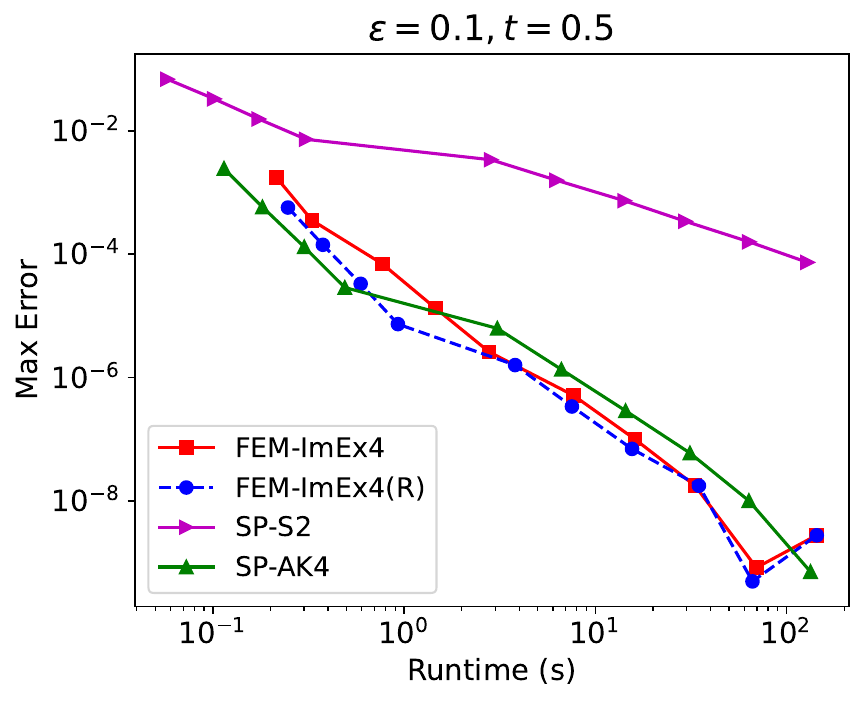}
            \caption{Cost measurement comparison.}
    \end{subfigure}
    \caption{Convergence rates and efficiency of different methods for the strong $\mathcal{O}(1)$ cubic focusing nonlinearity with $\epsilon = 0.2, 0.1$ and a nonzero initial phase. Meshing strategies for top row: $(\epsilon,\Delta x) = \left(0.2,\frac{1}{64}\right)$; bottom row: $(\epsilon,\Delta x) = \left(0.1,\frac{1}{128}\right)$.}
    \label{fig:Nonzero_phase_IC_cost_comp}
\end{figure}

Finally, we present convergence and work-precision diagrams for all of the fixed step size methods applied to the cases with $\epsilon=0.2$ and $\epsilon=0.1$.  We show results for the constant-phase initial data \eqref{constant-phase} in Figure \ref{fig:Zero_phase_IC_cost_comp}, and results for the varying-phase initial data \eqref{varying-phase}
in Figure \ref{fig:Nonzero_phase_IC_cost_comp}.
For these plots, we employ a finer spatial mesh compared to the mesh sizes used in the experiments above, and the spatial mesh is kept constant as the time step is refined. Specifically, for $\epsilon = 0.2, 0.1$, we choose $\Delta x = \frac{1}{64}, \frac{1}{128}$, respectively, for both initial data. The errors in the plots are computed as the maximum absolute error over the spatial grid, measured against a highly-accurate numerical reference solution. Convergence rates are estimated in the asymptotic regime and shown in the plots.  In each figure, the left panels show the solution convergence, and they verify that all the methods attain their expected order of accuracy. They also show a consistent improvement in accuracy when relaxation is applied.  The right panels are work-efficiency diagrams (error versus runtime), and the most efficient method corresponds to the curve that is closest to the bottom-left corner.  For $\epsilon=0.2$, the high-order splitting method is somewhat more efficient, while at $\epsilon=0.1$ all of the high-order methods yield similar efficiency.

\section{Preserving mass and energy through multiple relaxation}
\label{sec:multiple-relaxation-tests}
In the previous sections, relaxation was used to enforce conservation of mass only.
In this section, we use relaxation to conserve both mass and energy, and we study the behavior of the numerical solution and errors for long-time simulations.  To apply relaxation to both mass and energy,
we must use a semi-discretization that exactly conserves both quantities; we therefore use in this section
the FEM discretization described in Section \ref{Sec:Spatial_discretization}.
Here we first verify the semi-discrete conservation properties numerically.  Let $U(t)$ denote the exact solution
projected onto the finite element space.  We compute the time rate of change of the discrete mass or energy via
\begin{subequations}\label{Eq:time_deriv_invs}
    \begin{align}
    	\frac{d\eta_j(U(t))}{dt}  = \nabla \eta_j(U(t)) \cdot U'(t)   & = \nabla \eta_j(U(t)) \cdot f_{\textrm{FEM}} \;,  \\
    \end{align}   
\end{subequations}
where
\begin{align}
    \nabla \eta_1(U(t)) & = 2\Delta x \left[v_1(t),w_1(t),v_2(t),w_2(t), \cdots, v_m(t),w_m(t)\right]^T \;, 
\end{align}
and 
\begin{multline}
    \nabla \eta_2(U(t)) = \frac{2}{\Delta x} [(-v_{m-1}+2v_{1}-v_{2}),(-w_{m-1}+2w_{1}-w_{2}), (-v_{1}+2v_{2}-v_{3}), \\ 
    (-w_{1}+2w_{2}-w_{3}), \cdots, (-v_{m-2}+2v_{m-1}-v_{1}),(-w_{m-2}+2w_{m-1}-w_{1})
    ]^T  \\
    - 2 \beta \Delta x[(v_{1}^{2}+w_{1}^{2})v_{1},(v_{1}^{2}+w_{1}^{2})w_{1}, (v_{2}^{2}+w_{2}^{2})v_{2},(v_{2}^{2}+w_{2}^{2})w_{2}, \cdots, \\
    (v_{m-1}^{2}+w_{m-1}^{2})v_{m-1},(v_{m-1}^{2}+w_{m-1}^{2})w_{m-1}
    ]^T \;.
\end{multline}
For the 2-soliton and 3-soliton solutions, we take $12$ and $16$ grid points per unit length of the domain $[-35,35]$, respectively, and plot the time derivatives of two invariants using the exact solutions (see appendix~\ref{app:exact_multi_soliton_solution}) in \cref{Eq:time_deriv_invs} over the time domain $[0,T] = [0,20]$. 
\begin{figure}
 \centering
 \begin{subfigure}[b]{0.48\textwidth}
     \centering
     \includegraphics[width=\textwidth]{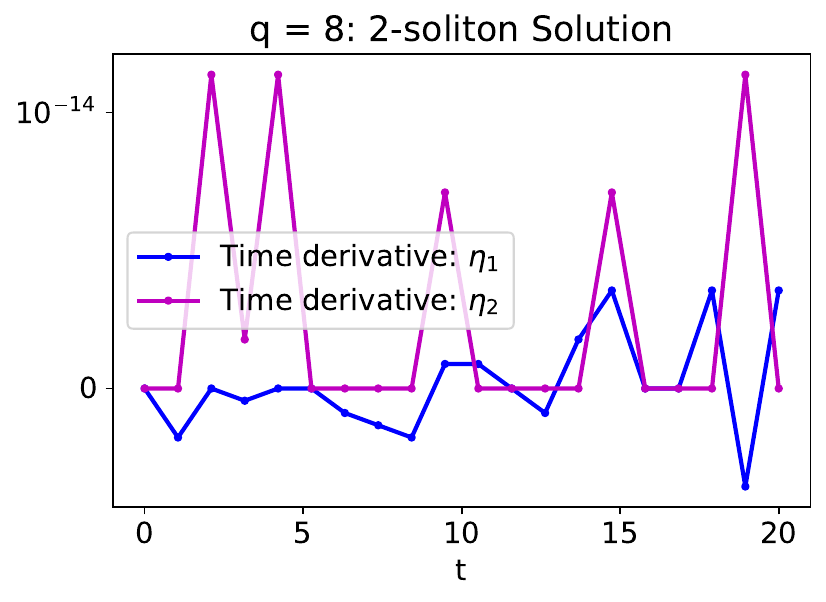}
 \end{subfigure}
 \hfill
 \begin{subfigure}[b]{0.48\textwidth}
     \centering
     \includegraphics[width=\textwidth]{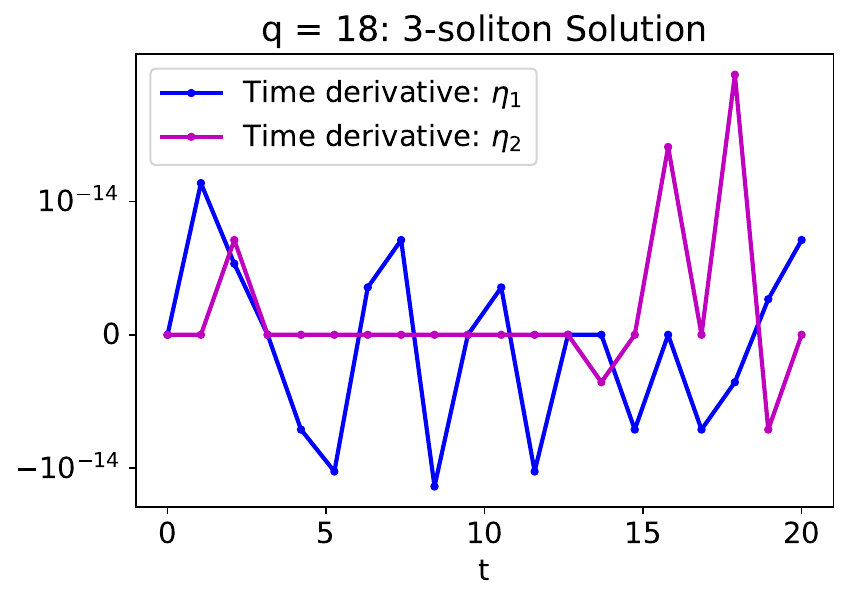}
 \end{subfigure}
	\caption{Time derivative of the discrete analogues of the invariants \eqref{Herbst_NLS_invariants} for the finite element semi-discretization given in \cref{ODE_system_FEM}.}
	\label{fig:time_deriv_inv}
\end{figure}
Figure ~\ref{fig:time_deriv_inv} confirms the conservation, up to rounding errors, of mass and energy by the finite element semi-discretization in \cref{ODE_system_FEM}.

In order to maintain the conservation of mass and energy in the fully discrete system, 
we apply relaxation to both quantities, as described in Section \ref{SubSubSec:MR_ImEx_Methods}.  The preservation of these quantities not only prevents the solution from growing arbitrarily, but it also provides a quantitative advantage in terms of error growth during long-term simulations. A comprehensive analysis of error propagation for the NLS equation with solitary waves was conducted in \cite{duran2000numerical}. The authors demonstrated that numerical methods that conserve the conserved quantities exhibit a better error propagation mechanism than nonconservative methods. The study considered the NLS equation with a relative equilibrium solution (i.e., a 1-soliton solution) and three nonlinear invariants, of which two are given in \eqref{Herbst_NLS_invariants}. The authors demonstrated analytically and experimentally that methods that conserve any two of these invariants lead to linear growth of errors over time, while non-conservative methods lead to quadratic error growth.  

    \begin{figure}
        \centering
        \includegraphics[width=\textwidth]{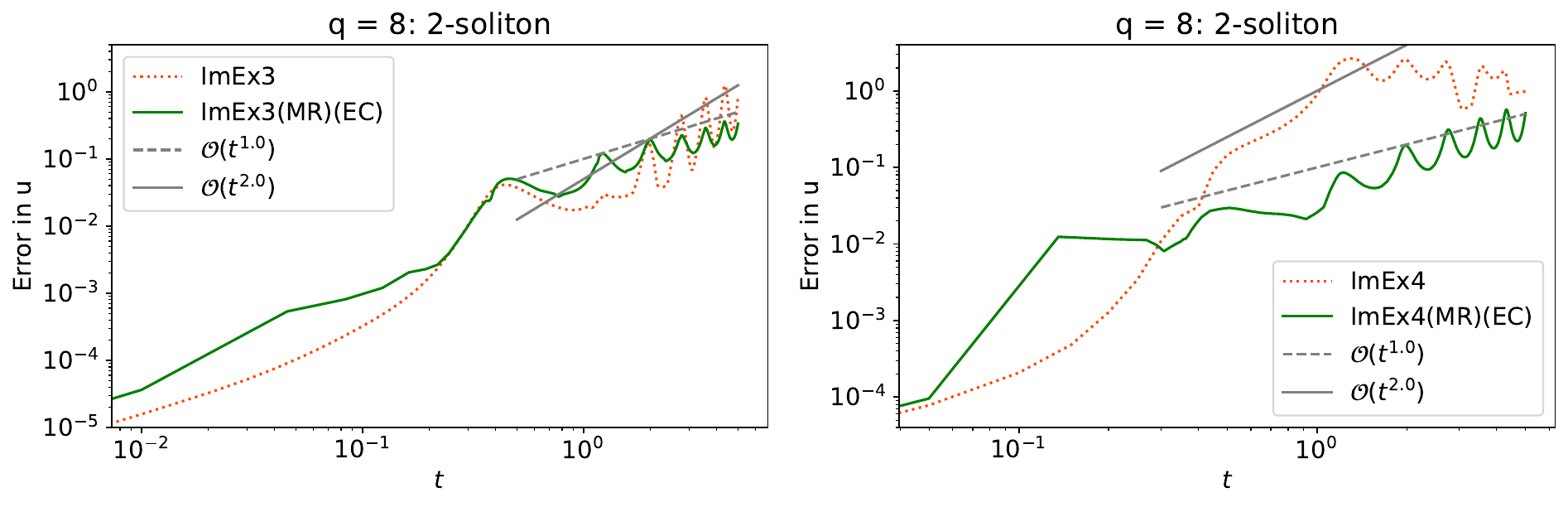}
        \caption{Error growth over time for a 2-soliton solution of the NLS \cref{Eq:Herbst_NLS}.  With FEM semi-discretization, relaxation is used to enforce conservation of $\eta_1$ and $\eta_2$.}
        \label{fig:q8_FEM_MRRK_SolErr_time}
    \end{figure}
    
    \begin{figure}
        \centering
        \includegraphics[width=\textwidth]{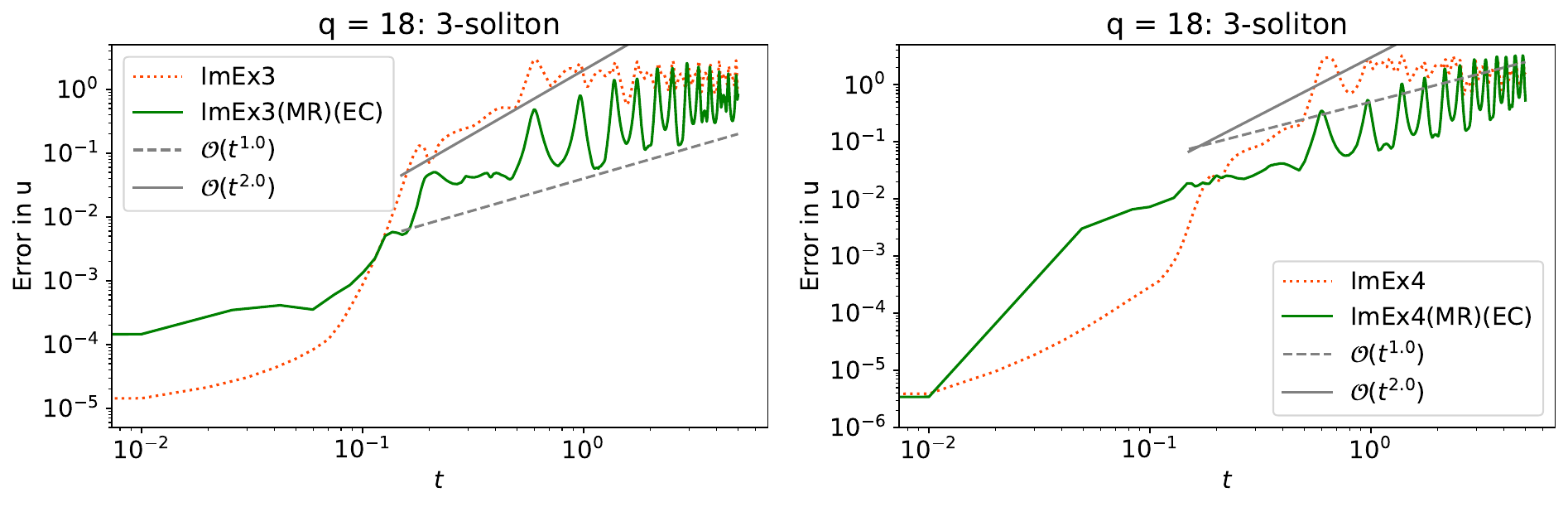}
        \caption{Error growth over time for a 3-soliton solution of the NLS \cref{Herbst_NLS_invariants}.  With FEM semi-discretization, relaxation is used to enforce conservation of $\eta_1$ and $\eta_2$.}
        \label{fig:q18_FEM_MRRK_SolErr_time}
    \end{figure}

It is natural to ask whether conservation of invariants will have a beneficial effect on more general (i.e., not just 1-soliton) solutions of the NLS equation, and a natural starting place is to consider solutions consisting of multiple solitons.  We seek to answer this question in this section.
To this end, we investigate the error behavior of the time-stepping strategy described above in terms of the invariants and solution error growth over time. Specifically, we apply two ImEx methods (ImEx3 and ImEx4) with and without the combination of relaxation and the hybrid time step size control.  We consider the 2- and 3-soliton solutions (see  Appendix~\ref{app:exact_multi_soliton_solution} for exact solutions) on the domain $[x_{L},x_{R}]=[-35,35]$ 
and integrate the semidiscretization \eqref{ODE_system_FEM} from $t=0$ to $t=5$. For the numerical simulation presented in this study, we set the error tolerance to $\tol = 10^{-4}$.

\begin{table}[ht]
    \centering
    \resizebox{12cm}{!}{%
    \begin{tabular}{ | c| l | c | c |}
    \hline
               & &        \multicolumn{2}{c|}{Maximum of}\\
   \hline
        Case & Methods & $|\eta_{1}\left(U(t))-\eta_{1}(U(0)\right)|$ & $|\eta_{2}\left(U(t))-\eta_{2}(U(0)\right)|$ \\
        \hlineB{2.5}
    \multirow{4}*{Two soliton} & FEM-ImEx3  & 9.50E-03 & 0.19 \\
        \cline{2-4}
                             & FEM-ImEx3(MR)(EC)  & 8.88E-16 & 7.11E-15 \\
        \cline{2-4}
                             & FEM-ImEx4  & 1.88E-01 & 2.64  \\
        \cline{2-4}
                             & FEM-ImEx4(MR)(EC)  & 1.11E-15 & 6.22E-15\\
    \hlineB{2.5} 
    \multirow{4}*{Three soliton} & FEM-ImEx3  & 8.57E-02 & 3.56  \\
        \cline{2-4}
                             & FEM-ImEx3(MR)(EC)  & 1.33E-15 & 1.78E-14  \\
        \cline{2-4}
                             & FEM-ImEx4  & 5.63E-02 & 2.81  \\
        \cline{2-4}
                             & FEM-ImEx4(MR)(EC)  & 1.55E-15 & 2.31E-14 \\
        \hline 
    \end{tabular}%
    }
    \caption{Maximum errors in invariants \eqref{Herbst_NLS_invariants} by different methods applied to a conservative finite element semi-discretization of the NLS equation with different n-solitons.}
    \label{Table:Invariant_Error}
\end{table}
\subsection{Bound State of Two Solitons}
 For the 2-soliton solution, we select a spatial grid consisting of $m=1120$ points and implement the ImEx3 and ImEx4 methods with time steps of $\Delta t = 0.01$ and $\Delta t = 0.05$, respectively. To incorporate the multiple relaxation approach in conjunction with the hybrid step size control for the ImEx methods, we initialize the time-stepping with the same time steps mentioned above. The maximum error in the conserved quantities for the 2-soliton solutions is presented in Table~\ref{Table:Invariant_Error}, indicating that all the modified methods preserve both invariants. 
 
 Figure~\ref{fig:q8_FEM_MRRK_SolErr_time} presents the solution error growth for all methods in the case of 2-soliton solutions. 
 In these plots there is a periodic variation in the error due to the time-periodic terms in the exact solution.
 While the improvement in solution accuracy is particularly pronounced for the 4th-order methods, it can be seen
in all cases that the conservative (multiple-relaxation) methods exhibit linear error growth, while the baseline methods yield errors that increase approximately quadratically over time.  In all of these tests, the error will eventually \emph{saturate}, or become so large that it stops growing, which can be seen in some of the plots at late times.
\subsection{Bound State of Three Solitons}
To capture the sharp gradients in the 3-soliton solution, we use an even finer spatial grid with $m = 2240$ points and use an initial time step size of $\Delta t = 0.01$ for all the methods. The maximum error in the conserved quantities for the 3-soliton solutions is provided in Table~\ref{Table:Invariant_Error}, confirming that both invariants are conserved by the modified ImEx methods. Additionally, Figure~\ref{fig:q18_FEM_MRRK_SolErr_time} illustrates the solution error growth of all methods for the case of 3-soliton solutions. The modified ImEx methods with step size control that conserve the two invariants exhibit a slighly superlinear error growth rate, while the baseline methods have a rate with exponent slightly larger than $2$ for long times.  Conservation again yields a significant improvement in the error growth rate, and correspondingly
in the size of the error itself for long times, although the errors again eventually saturate.

\begin{figure}[ht]
 \centering
 \begin{subfigure}[t]{0.48\textwidth}
     \centering
     \includegraphics[width=\textwidth]{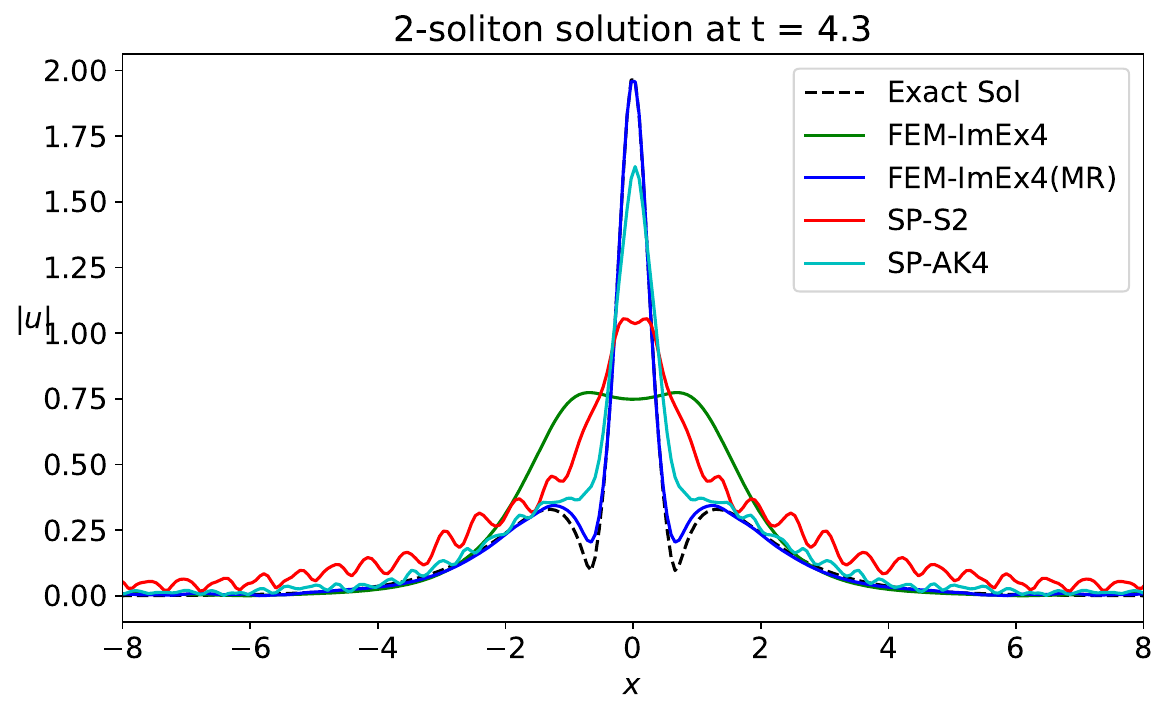}
     \caption{Two-soliton solution with $\Delta x = 1/16$ and $\Delta t = 1/20$.}
 \end{subfigure}
 \hfill
 \begin{subfigure}[t]{0.48\textwidth}
     \centering
     \includegraphics[width=\textwidth]{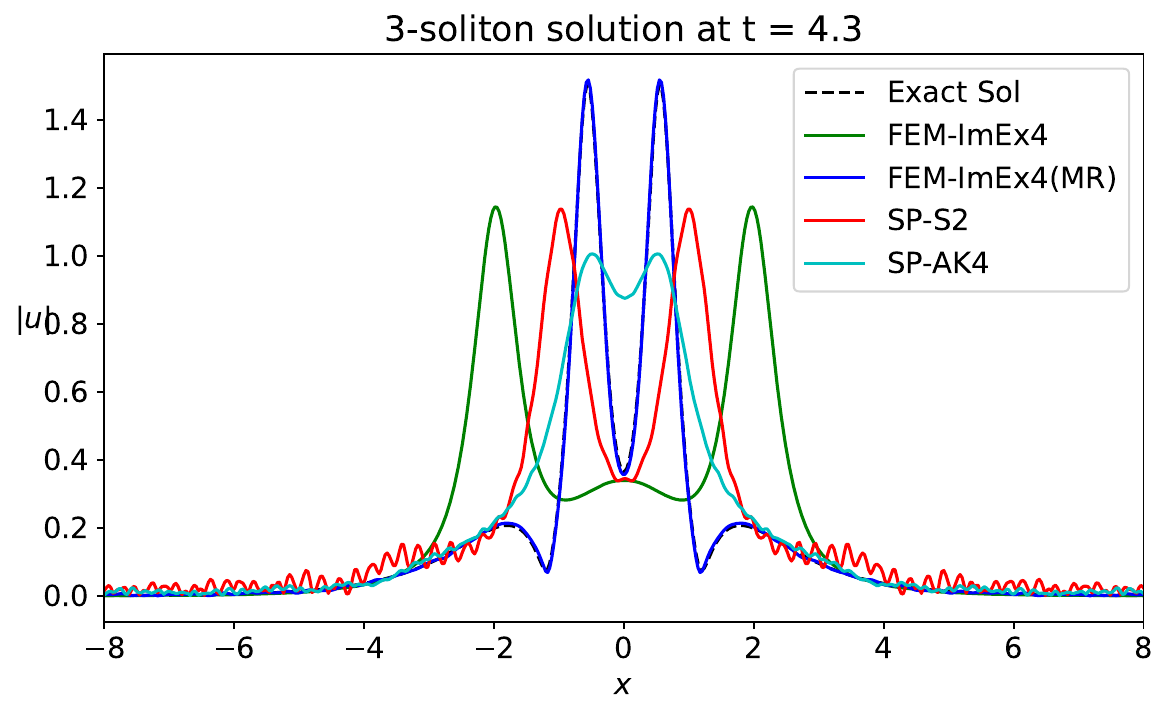}
     \caption{Three-soliton solution with $\Delta x = 1/32$ and $\Delta t = 1/100$.  Note that the exact solution
              is nearly indistinguishable from that of FEM-ImEx4(MR).}
 \end{subfigure}
	\caption{Comparison of numerical approximations for the multi-soliton solutions of the NLS \cref{Eq:Herbst_NLS} using different methods. The figure compares time-splitting pseudo-spectral methods with a full discretization that employs finite element approximation in space and ImEx methods with and without relaxation in time.}
	\label{fig:Sol_Comparison_Herbst_NLS}
\end{figure}

To give a clearer visual comparison of the improvement in the error by conservative methods, in Figure~\ref{fig:Sol_Comparison_Herbst_NLS} we plot solutions in space 
obtained using each of the time-stepping methods considered in this work.  Here we have used the same fixed step size with all methods, in order to show
clearly that the differences are related to conservation properties and not to adaptive time stepping.  For these comparisons we have deliberately chosen
 a late time at which the solution exhibits large values of the spatial derivatives (what almost appear to be kinks in the exact solution).

These figures dramatically illustrate the improved accuracy obtained with the conservative methods.  It should be noted that the error for all methods is
roughly time periodic (due to the nature of the exact solution), and the difference between numerical solutions is less pronounced at some times.
Furthermore, it should be noted that the FEM semi-discretization has a higher computational cost compared to the pseudospectral discretization, since
in the former the solution of a linear algebraic system is required at each Runge-Kutta stage.  At the same time, it is remarkable that much better
accuracy is obtained with a 2nd-order FEM semi-discretization compared to the spectral semi-discretization, and this highlights even further the
importance of conservation. The inadequate performance of the splitting methods can be explained by following the work in \cite{verwer1986conservative}. Essentially, these methods treat the dispersive and nonlinear components separately rather than concurrently, resulting in an imbalance between dispersive and nonlinear forces. The dispersive characteristics of the linear component cause the solution to disperse when solved using a fractional step. However, the subsequent fractional step for solving the nonlinear component does not adequately mitigate this dispersion phenomenon since there is no coupling between adjacent grid points.

\section{Conclusion and Perspectives}
\label{Sec:conclusion}
Pseudospectral discretization combined with low-order operator splitting has long been the standard
for numerical solution of NLS and related equations.  However, these methods are known to be
inefficient in some regimes \cite{bao2003numerical}, and recent work has highlighted the potential advantages
of higher-order time stepping \cite{antoine2016high}.  Here we have provided an essentially explicit approach
that has the additional advantages of both mass conservation and high-order accuracy.
We have shown that this combination results in a generally competitive method that has substantial
advantages in the semiclassical limit of small $\epsilon$.  When combined with adaptive time stepping,
it can yield huge improvements in computational efficiency and accuracy.

In a more theoretical vein, we have explored the effects of numerical conservation
on long-time error growth.  The landmark result in this area applies only to the
rather simple 1-soliton solution, so we have studied experimentally the case of
multiple solitons, finding that a similar qualitative improvement (from quadratic
to linear error growth in time) is observed.

Future work may explore the application of our proposed time-stepping approach to other types of nonlinear Schrödinger equations, such as the multi-component nonlinear Schrödinger or Gross–Pitaevskii equation, which is used for modeling Bose-Einstein condensation.  Additionally, it seems of interest to study
numerical error growth for multi-soliton and other solutions, to see if the
observations here can be supported through analysis. 
\appendix
\section{Exact Multi-soliton Solutions}
\label{app:exact_multi_soliton_solution}
    Following the work in \cite{miles1981envelope}, explicit formulas for the exact solutions corresponding to 1-soliton, 2-soliton, and 3-soliton solutions for the nonlinear Schr\"odinger equation ~\eqref{Eq:Herbst_NLS}  can be written down, respectively, as follows: 

   \begin{align}
        u_1(x,t) = \frac{2 e^{i t} }{e^{x}+e^{-x}}\;,
    \end{align}
    
    \begin{align}
        u_2(x,t) = \frac{2 e^{x+9 i t} \left(3 e^{2 x+8 i t}+3 e^{4 x+8 i t}+e^{6 x}+1\right)}{3 e^{4 (x+4 i t)}+4 e^{2 x+8 i t}+4 e^{6 x+8 i t}+e^{8 x+8 i t}+e^{8 i t}+3 e^{4 x}}\;, \ \text{and}
    \end{align}
    
     \begin{align}
        u_3(x,t)
            & = \frac{
                    \begin{aligned}
                        80 e^{7 (x+7 i t)}+2e^{x+25 i t}+16 e^{3 x+33 i t}+36 e^{5 x+33 i t}+20 e^{5 x+49 i t} + \\
                        32 e^{7 x+25 i t}+10 e^{9 x+9 i t}+90 e^{9 x+41 i t}+40 e^{9 x+57 i t}+32 e^{11 x+25 i t} + \\
                        80 e^{11 x+49 i t}+32 e^{13 x+33 i t}+20 e^{13 x+49 i t}+16 e^{15 x+33 i t}+2e^{17 x+25 i t} 
                    \end{aligned}
                    }
                    {
                    \begin{aligned}
                        64 e^{12 (x+2 i t)}+36 e^{8 (x+3 i t)}+18 e^{4 (x+4 i t)}+64 e^{6 (x+4 i t)}+45 e^{10 (x+4 i t)}+10 e^{12 (x+4 i t)}+ \\
                        45 e^{8 (x+5 i t)}+18 e^{4 (x+8 i t)}+10 e^{6 (x+8 i t)}+9 e^{2 x+24 i t}+45 e^{8 x+8 i t}+45 e^{10 x+8 i t}+36 e^{10 x+24 i t}+ \\
                        18 e^{14 x+16 i t}+18 e^{14 x+32 i t}+9 e^{16 x+24 i t}+e^{18 x+24 i t}+e^{24 i t}+10 e^{6 x}+10 e^{12 x}
                    \end{aligned}
                    }
                    \;.
    \end{align}

\vspace{1.5em}
\bibliographystyle{plain}
\bibliography{refs}

\vspace{1.5em}
\end{document}